\documentclass[a4paper]{article}
\usepackage[utf8]{inputenc}
\usepackage[ngerman, english]{babel}
\usepackage{amsmath,amsthm,amssymb,amsfonts}
\usepackage{mathrsfs}
\usepackage{ifthen}
\usepackage{enumerate}
\usepackage{hyphsubst}
\usepackage{xcolor}
\usepackage{graphicx}
\usepackage{psfrag}
\usepackage{pst-all}

\newtheorem{theorem}{Theorem}[section]
\newtheorem{lemma}[theorem]{Lemma}

\newtheorem{proposition}[theorem]{Proposition}
\newtheorem{definition}[theorem]{Definition}

\newtheorem{remark}[theorem]{Remark}

\newcommand{\nv}{\boldnu}
\newcommand{\domain}{D}
\newcommand{\angvel}{\Omega}

\newcommand{\dens}{\rho}
\newcommand{\densl}{\ul{\rho}}

\newcommand{\pres}{p}
\newcommand{\grav}{\phi}
\newcommand{\dislag}{\boldu}
\newcommand{\dislagscal}{u}
\newcommand{\graveu}{\psi}
\newcommand{\bflow}{\boldb}
\newcommand{\qcoeff}{\boldq}
\newcommand{\sounds}{c_s}
\newcommand{\soundsl}{\ul{c_s}}

\newcommand{\damp}{\gamma}
\newcommand{\dampl}{\ul{\gamma}}

\newcommand{\ses}{a}
\newcommand{\sescow}{\ses_\mathrm{Cow}}

\newcommand{\op}{A}
\newcommand{\opAcow}{\op_\mathrm{Cow}}

\newcommand{\hs}{\boldX}
\newcommand{\hsV}{\boldV}
\newcommand{\hsW}{\boldW}
\newcommand{\hsZ}{\boldZ}
\newcommand{\hsgeneric}{Y}
\newcommand{\ugeneric}{y}
\newcommand{\hsgrav}{\tilde H^1_*}

\newcommand{\blo}{L}
\newcommand{\compv}{\boldv}
\newcommand{\compvp}{v_0}
\newcommand{\compw}{\boldw}

\newcommand{\compz}{\boldz}
\newcommand{\PV}{P_\hsV}
\newcommand{\PW}{P_\hsW}
\newcommand{\PZ}{P_\hsZ}
\newcommand{\PHi}{P_{H_1}}
\newcommand{\PHii}{P_{H_2}}

\newcommand{\mati}{\ul{\ul{m_1}}}
\newcommand{\matii}{\ul{\ul{m_2}}}
\newcommand{\Xint}{\hs_\mathrm{int}}
\newcommand{\hsVint}{\boldV_\mathrm{int}}
\newcommand{\hsWint}{\boldW_\mathrm{int}}
\newcommand{\hsZint}{\boldZ_\mathrm{int}}
\newcommand{\PVint}{P_{\hsVint}}
\newcommand{\PWint}{P_{\hsWint}}
\newcommand{\PZint}{P_{\hsZint}}
\newcommand{\Vext}{V_\mathrm{ext}}
\newcommand{\sescp}{\ses_\mathrm{cp}}
\newcommand{\Acp}{A_\mathrm{cp}}
\newcommand{\Tcp}{T_\mathrm{cp}}
\newcommand{\Tint}{T_\mathrm{int}}
\newcommand{\Kf}{K_3}

\newcommand{\dd}{\mathrm{d}}
\newcommand{\ol}[1]{\overline{#1}}
\newcommand{\ul}[1]{\underline{#1}}
\newcommand{\spl}{\langle}
\newcommand{\spr}{\rangle}
\newcommand{\bpm}{\begin{pmatrix}}
\newcommand{\epm}{\end{pmatrix}}

\renewcommand{\div}{\operatorname{div}}

\DeclareMathOperator{\curl}{curl}
\DeclareMathOperator{\hess}{Hess}

\DeclareMathOperator{\supp}{supp}

\DeclareMathOperator{\ran}{ran}

\DeclareMathOperator{\sign}{sgn}

\DeclareMathOperator{\numran}{num ran}

\newcommand{\boldnu}{\boldsymbol{\nu}}
\newcommand{\boldxi}{\boldsymbol{\xi}}

\newcommand{\setC}{\mathbb{C}}

\newcommand{\setR}{\mathbb{R}}

\newcommand{\boldb}{\mathbf{b}}

\newcommand{\boldf}{\mathbf{f}}
\newcommand{\boldg}{\mathbf{g}}

\newcommand{\boldq}{\mathbf{q}}

\newcommand{\boldu}{\mathbf{u}}
\newcommand{\boldv}{\mathbf{v}}
\newcommand{\boldw}{\mathbf{w}}
\newcommand{\boldx}{\mathbf{x}}

\newcommand{\boldz}{\mathbf{z}}

\newcommand{\boldG}{\mathbf{G}}
\newcommand{\boldH}{\mathbf{H}}

\newcommand{\boldL}{\mathbf{L}}

\newcommand{\boldV}{\mathbf{V}}
\newcommand{\boldW}{\mathbf{W}}
\newcommand{\boldX}{\mathbf{X}}

\newcommand{\boldZ}{\mathbf{Z}}

\newcommand{\calT}{\mathcal{T}}

\title{On the treatment of exterior domains for the time-harmonic equations of stellar oscillations\footnote{This work was partially funded by the Deutsche Forschungsgemeinschaft (DFG, German Research Foundation) – Project-ID 432680300 - SFB 1456}}
\author{Martin Halla\footnote{Max-Planck-Institut f\"ur Sonnensystemforschung, Justus-von-Liebig-Weg 3, 37077 G\"ottingen, Deutschland (halla@mps.mpg.de)}$~^,$
\footnote{Institut f\"ur Numerische und Angewandte Mathematik, Georg-August Universit\"at G\"ottingen, Lotzestra\ss e 16-18, 37083 Göttingen, Deutschland}}
\date{\today}

\begin{document}
\maketitle
\begin{abstract}
\noindent
In a recent article we started to analyze the time-harmonic equations of stellar oscillations.
As a first step we considered bounded domains together with an essential boundary condition and established the well-posedness of the equation.
In this article we consider the physical relevant case of the domain being $\setR^3$.
We discuss the treatment of the exterior domain, and show how to couple the two parts to obtain a well-posedness result.
Further, for the Cowling approximation (which neglects the Eulerian perturbation of gravity) we derive a scalar equation in the atmosphere, couple it to the vectorial interior equation, and prove the well-posedness of the new system.
This coupled system has the big advantages that it simplifies the construction of approximating transparent boundary conditions and leads to significant less degrees of freedom for discretizations.\\


\end{abstract}

\section{Introduction}\label{sec:introduction}

In this article we study the time-harmonic linear equations of stellar oscillations \cite{HallaHohage:21, LyndenBOstriker:67} (with the phase convention $e^{-i\omega t}$)
\begin{subequations}\label{eq:PDE}
\begin{align}\label{eq:PDE-dislag}
\begin{split}
-\nabla \big(\sounds^2 \dens\div \dislag
+\nabla \pres\cdot\dislag\big)
+\nabla \pres \div \dislag
+\hess(\pres)\dislag
& \phantom{= \dens\boldf \quad\text{in }\setR^3,}\\
-\dens\hess(\grav)\dislag
-\dens\nabla\graveu
- \dens\big(\omega +i\partial_{\bflow}+i\angvel\times \big)^2 \dislag
-i\omega\damp\dens\dislag
&= \dens\boldf \quad\text{in }\setR^3,
\end{split}\\
\label{eq:PDE-gravlag}
-\frac{1}{4\pi G}\Delta \graveu +\div(\dens\dislag) &= \phantom{\dens}0 \quad\text{in }\setR^3.
\end{align}
together with the decay conditions
\begin{align}\label{eq:decay}
\dislag \in H(\div;\setR^3)
\qquad\text{and}\qquad
\nabla\graveu \in \boldL^2(\setR^3)
\end{align}
\end{subequations}
and a Gauge condition for $\graveu$ to cancel out the constant functions.
Thereat $\dislag$ is the Lagrangian perturbation of displacement and $\graveu$ is the (scaled) Eulerian perturbation of the 
gravitational background potential $\phi$.
The equation is formulated in a frame rotating at constant angular velocity 
$\angvel\in\mathbb{R}^3$ with the star.
Further, $\dens, \pres, \sounds, \bflow, G$ and $\boldf$ denote density, pressure, sound speed, background velocity, gravitational constant and sources, 
$\partial_\bflow := \sum_{l=1}^3 \bflow_l\partial_{\boldx_l}$  
denotes the directional derivative in direction $\bflow$, 
$\hess(v)$ the Hessian matrix of a scalar function $v=\pres,\grav$, and damping is modeled by the term $-i \omega \damp \dens \dislag$ 
with damping coefficient $\damp$.
Equation~\eqref{eq:PDE} with $\damp\equiv0$ was first derived in \cite{LyndenBOstriker:67} and appears as Lagrangian linearization of the time-dependent non-linear Euler equations around a stationary solution $(\bflow,\dens,\pres,\grav)$, i.e.\ $(\bflow,\dens,\pres,\grav)$ satisfy
\begin{subequations}\label{eq:equi}
\begin{align}\label{eq:motion-equi}
\dens\big(\partial_\bflow\bflow +2\angvel\times\bflow+\angvel\times(\angvel\times\boldx)
-\nabla\grav\big)+\nabla \pres &= 0 \quad\text{in }\setR^3,\\
\label{eq:cont-equi}
\div(\dens\bflow)&=0 \quad\text{in }\setR^3,\\
\label{eq:gravity-equi}
-\Delta \grav -4\pi G\dens &= 0 \quad\text{in }\setR^3.
\end{align}
\end{subequations}
The linearized equations can be reduced to a system for $(\dislag,\graveu)$, and
subsequently the term $\damp \dens(- i \omega) \dislag$ is included into \eqref{eq:PDE-dislag} to model damping effects, which are believed to be caused mainly by radiative damping and interaction with turbulent convection (see \cite{UOAS:89}).
This particular choice of the damping term was first proposed in \cite{HallaHohage:21}, because different to other reasonable and simple damping models it stabilizes the equation.
For more refined damping models which involve non-local terms we refer e.g.\ to \cite{KK:64}.
Note that \eqref{eq:PDE-dislag} without damping, rotational and gravitational terms ($\damp\equiv0$, $\angvel=0$, $\grav\equiv\graveu\equiv0$) was first derived by Galbrun \cite{Galbrun:31} and is conveniently referred to as Galbrun's equation.
In this form the equation is used in aeroacoustics to model and eventually reduce noise caused by moving objects such as aircraft engines (see, e.g.\ \cite{GLT:04}).

In a preceding article \cite{HallaHohage:21} we established the well-posedness of \eqref{eq:PDE} in a bounded interior domain $D_\mathrm{int}$ with the boundary condition $\nv\cdot\dislag=0$ on $\partial D_\mathrm{int}$ (see also \cite{BonnetBDMercierMillotPernetPeynaud:12,HaeggBerggren:19} for results in this direction).
The present paper is devoted to the treatment of the exterior domain $D_\mathrm{Atmo}:=\setR^3\setminus \ol{D_\mathrm{int}}$.
Since \eqref{eq:PDE-dislag} is formulated in a rotating frame the rotation of which is aligned to the rotation of the star, it makes sense physically to assume that the background flow vanishes in the atmosphere (i.e.\ $\bflow=0$ in $D_\mathrm{Atmo}$), and we will embrace this assumption in the entire article.
In general, solutions to wave equations in open domains do not decay fast enough and it is necessary to apply special mathematical tools in the exterior domain to obtain Fredholmness results.
We refer to \cite{Kim:14,Halla:20PML} for complex scaling/perfectly matched layer methods, to \cite{HohageNannen:09,Halla:16} for Hardy space/pole condition methods, and to \cite{Hagstrom:99,Givoli:92} for radiation boundary conditions.
We mention \cite{BecacheBonnetBDLegendre:06,RetkaMarburg:13,BaccoucheTaharMoreau:16} for complex scaling and infinite element methods for Galbrun's equation with uniform flows.
For spherical symmetric backgrounds and idealistic stellar parameters a modal analysis for a simplified Galbrun's equation ($\angvel=0$, $\bflow\equiv0$, $\graveu\equiv0$) is reported in \cite{BarucqFaucherFournierGizonPham:21} by means of the Liouville transform.
In our stellar context the equation contains a damping term and thus it seems reasonable that an intricate radiation condition can be avoided and be replaced by a simple decay condition.
However, at this point this notion is purely intuitive and needs to be justified properly.
The damping and the locality of the background flow for our stellar configuration simplify the analysis of \eqref{eq:PDE-dislag} in the atmosphere.
Nevertheless, different to equations which model wave propagation e.g.\ in the air or earth crust, in our stellar context we have to deal with non-homogeneous parameters in the exterior domain.
In particular, in the atmosphere of the star the density $\dens$ decays to zero for increasing radius.

A key ingredient of our analysis is that the part of the sesquilinearform associated to \eqref{eq:PDE} which involves only $\dislag$-components can in the atmosphere be formulated as
\begin{align*}
\spl \sounds^2\dens (\div+\qcoeff\cdot)\dislag, (\div+\qcoeff\cdot)\dislag' \spr_{L^2(D_\mathrm{Atmo})}
- \spl \dens(i\omega\damp+\matii)\dislag, \dislag' \spr_{\boldL^2(D_\mathrm{Atmo})}
\end{align*}
with a vectorial parameter $\qcoeff$ and a selfadjoint matrix function $\matii$.
From this representation it follows that this part of the sesquilinear form is coercive.
Thus the main work of our analysis is to develop a technique to couple this coercivity result in the atmosphere to the ideas developed in \cite{HallaHohage:21} for the interior part.
To achieve this we will introduce a transition layer between the interior domain and the atmosphere.
In addition, from the coercivity in the atmosphere it follows that indeed the decay condition is sufficient to describe the behavior of outgoing solutions at infinity.
We report our first main result, the well-posedness of \eqref{eq:PDE} in Theorem~\ref{thm:bij-full}.

An additional main achievement of this article concerns the Cowling approximation (which neglects the Eulerian perturbation of gravity $\graveu$) with spherical symmetric parameters in the atmosphere.
In the atmosphere we obtain a representation of $\dislag$ in terms of a \emph{scalar} potential \eqref{eq:UeqGrad}.
Thereupon, we derive an equation for the potential, formulate a system of equations which couples the interior and the exterior part \eqref{eq:PDE-cp}, and prove its well-posedness in Theorem~\ref{thm:wTc-coupled} and Proposition~\ref{prop:bij-cp}.
This new system achieves a significant simplification for the construction of numerical transparent boundary condition methods, and a compelling reduction of degrees of freedom in the atmosphere.
Recall that for (finite element) discretizations it is convenient to truncate the open domain to a bounded one and to impose an approximated transparent boundary condition at the artificial boundary.
The simplest possibility is to choose a homogeneous essential boundary condition, and from the coercivity in the atmosphere one can even deduce that for increasing domain sizes this approximation converges (see \cite{Halla:20PML,Halla:19Diss} for such analysis techniques).
However, the convergence speed depends on the physical damping parameter $\damp$, which cannot be tuned, and thus this approach is too expensive for practical applications.
Up to now numerical transparent boundary conditions for stellar equations have only been reported for simplified scalar equations and we refer to \cite{BarucqFaucherPham:20} for radiation boundary conditions and to \cite{HohageLehrenfeldPreuss:20} for learned infinite elements.

The remainder of this article is structured as follows.
In Section~\ref{sec:preliminaries} we formulate our assumptions on the physical parameters, we derive the variational formulation of \eqref{eq:PDE}, we set our notation and recall some common definitions.
In Section~\ref{sec:cowling} we consider the so-called Cowling approximation of \eqref{eq:PDE}, which neglects the Eulerian perturbation of gravity $\graveu$ in $\eqref{eq:PDE-dislag}$ and is solely an equation for the Lagrangian displacements $\dislag$.
The Cowling approximation already contains most of the mathematical difficulties to study \eqref{eq:PDE}.
In Theorem~\ref{thm:decomposition} we report a topological decomposition of the Hilbert space similar to \cite[Theorem~3.5]{HallaHohage:21}, and subsequently we report the weak T-coercivity and the bijective of the operator under investigation in Theorem~\ref{thm:wTc-cow} and Proposition~\ref{prop:bij-cow}.
The main work of this section is to construct a suitable transition for the T-operator between the interior domain and the atmosphere.
In Section~\ref{sec:full} we extend our results to the full equation \eqref{eq:PDE} in Theorem~\ref{thm:bij-full}.
Since the off-diagonal operators which couple the equations for $\dislag$ and $\graveu$ involve integrals over $\setR^3$, these operators are not compact, and hence the analysis requires some new ideas compared to \cite{HallaHohage:21}.
In Section~\ref{sec:cp} we consider the Cowling approximation and spherical symmetric parameters in the atmosphere.
In the atmosphere we derive a scalar equation for a potential of $\dislag$, formulate a coupled system and prove its weak T-coercivity and bijectivity in Theorem~\ref{thm:wTc-coupled} and Proposition~\ref{prop:bij-cp}.
We close this article with a conclusion and outlook in Section~\ref{sec:conclusion}.

\section{Preliminaries}\label{sec:preliminaries}

In Section~\ref{subsec:assumptions} we formulate our assumptions on the physical parameters,
in Section~\ref{subsec:varform} we derive the variational formulation of \eqref{eq:PDE},
and in Section~\ref{subsec:commondefs} we recall some common definitions.

\subsection{Basic assumptions}\label{subsec:assumptions}
We denote the spatial coordinate as $\boldx\in\setR^3$.
For $r>0$ let $B_r:=\{\boldx\in\setR^3\colon |\boldx|<r\}$.
For $r_2>r_1>0$ let $A_{r_1,r_2}:=B_{r_2}\setminus\ol{B_{r_1}}$.
Let $\omega\in\setR$ and $\angvel\in\setR^3$.
Let $\dens\in L^\infty(\setR^3,\setR)$ be such that
\begin{align}\label{eq:densl}
\densl_r:=\inf_{\boldx\in B_r}\dens(\boldx) > 0
\end{align}
for each $r>0$ and $\sounds,\damp\in L^\infty(\setR^3,\setR)$ be such that
\begin{align}\label{eq:box_constraints}
\begin{split}
\soundsl:=\inf_{\boldx\in\setR^3}\sounds(\boldx)>0
\qquad\text{and}\qquad
\dampl:=\inf_{\boldx\in\setR^3}\damp(\boldx)>0.
\end{split}
\end{align}
In particular, for stellar models we have to consider that the density tends to zero for increasing radius: $\lim_{r\to+\infty}\inf_{\boldx\in B_r^c}\dens(\boldx)=0$.
For example, the standard model S of \cite{C-DEtal:96} for the sun assumes in the atmosphere $\dens(\boldx)=Ce^{-\alpha|\boldx|}$ with positive constants $C,\alpha$, and we refer to \cite{VAL:81} for alternative models.
Further, let $\pres, \grav \in W_\mathrm{loc}^{2,\infty}(\setR^3,\setR)$ and
\begin{align}\label{eq:q-mati}
\qcoeff:=\frac{1}{\sounds^2\dens}\nabla\pres \qquad\text{and}\qquad
\mati:=-\dens^{-1}(\hess(\pres)-\dens\hess(\grav)-\sounds^2\dens\,\qcoeff\qcoeff^\top).
\end{align}
We assume that $\qcoeff\in L^\infty(\setR^3,\setR^3)$ and $\mati\in L^\infty(\setR^3,\setR^{3\times3})$.
For later use we also define the matrix function
\begin{align}\label{eq:matii}
\matii:=\mati+(\omega+i\angvel\times)^*(\omega+i\angvel\times).
\end{align}
For a Lipschitz domain $D\subset\setR^3$ we introduce the weighted spaces $L^2_\dens(D,\setC)$ and $L^2_\dens(D,\setC^3)$ with scalar products
\begin{align*}
\spl \dislagscal,\dislagscal' \spr_D:=\spl \dislagscal,\dislagscal' \spr_{L^2_\dens(D)}
:=\int_{D} \dens\dislagscal\ol{\dislagscal'} \,\dd\boldx, \quad
\spl \dislag,\dislag' \spr_D:=\spl \dislag,\dislag' \spr_{(L^2_\dens(D))^3}:=\int_{\setR^3} \dens\dislag\cdot\ol{\dislag'} \,\dd\boldx,
\end{align*}
for scalar functions $\dislagscal,\dislagscal'$ and vectorial functions $\dislag,\dislag'$, whereby $\ol{\cdot}$ denotes the complex conjugation. Since we use the same symbol for the scalar products of scalar and vectorial functions the notation is overloaded, but its meaning will always be clear from the arguments.
Further, for $D=\setR^3$ we set $\spl\cdot,\cdot\spr:=\spl\cdot,\cdot\spr_{\setR^3}$.
Henceforth we consider $\setR^3$ as the default domain for all function spaces and suppress the dependency in the notation, if the domain equals $\setR^3$.
Thus we write e.g.\ $L^2=L^2(\setR^3)$ and so on.
If we do not explicitly indicate a particular field, all spaces are over $\setC$, e.g.\ $L^2=L^2(\setR^3)=L^2(\setR^3;\setC)$.
Further, for any space $\hsgeneric$ of scalar functions $\ugeneric\colon\domain\to\setC$ we set $\mathbf{\hsgeneric}:=(\hsgeneric)^3$.
We denote the three-by-three identity matrix as $I_{3\times3}$.
For a scalar or matrix function $\sigma$ we denote the multiplication operator with symbol $\sigma$ as $M_\sigma$.
We denote the directional derivative in direction $\bflow\in\setR^3$ as
\begin{align*}
\partial_\bflow := \sum_{l=1}^3 \bflow_l\partial_{\boldx_l} = \bflow\cdot\nabla.
\end{align*}
Next we formulate assumptions on the flow $\bflow\in L^{\infty}(\setR^3,\setR^3)$ and recall some results from \cite{HallaHohage:21}.
Let $\div(\dens\bflow)\in L^2(\setR^3,\setR)$.
Thence we are able to well define the weak derivative
$\dens\partial_\bflow$ through
\begin{align*}
\spl \dens\partial_\bflow \dislag,\dislag' \spr :=
-\spl \dens\dislag,\partial_\bflow\dislag' \spr
-\spl\div(\dens\bflow)\dislag,\dislag'\spr
\end{align*}
for $\dislag'\in C^\infty_0(\setR^3,\setC^3)$ and so we set $\partial_\bflow:=\dens^{-1}(\dens\partial_\bflow)$.
This way we can define Sobolev spaces like $Y=\{\dislag\in \boldL^2_\dens\colon\partial_\bflow\dislag\in \boldL^2_\dens\}$, $\spl \boldu,\boldu'\spr_Y=\spl \boldu,\boldu'\spr+\spl \partial_\bflow \boldu,\partial_\bflow \boldu'\spr$ and their completeness follows as in \cite[Lemma~2.1]{HallaHohage:21}.
Since Equation~\eqref{eq:PDE-dislag} is formulated in a frame which is rotating together with the star at angular velocity $\Omega$, it is reasonable to assume that the flow $\bflow$ is local and vanishes in the atmosphere.
Hence, we assume that there exists a constant $r_1>0$ such that
\begin{align}
\supp\bflow:=\ol{\{\boldx\in\setR^3\colon\bflow(\boldx)\neq0\}}\subset B_{r_1}.
\end{align}
We note the general integration by parts formula
\begin{align*}
\spl \dens \partial_\bflow \dislag,\dislag' \spr_{\boldL^2(D)}
&=-\spl \dens \dislag, \partial_{\ol{\bflow}} \dislag' \spr_{\boldL^2(D)}
-\spl \div(\dens\bflow) \dislag, \dislag' \spr_{\boldL^2(D)}
+ \int_{\partial\domain} (\nv\cdot\bflow) (\dislag \cdot \ol{\dislag'}) \,\dd S
\end{align*}
for a Lipschitz domain $D\subset\setR^3$, and hence with our assumptions on $\bflow$ and \eqref{eq:cont-equi} it follows
\begin{align}\label{eq:IntbyPartsBflow}
\spl i\partial_\bflow \dislag,\dislag' \spr
&=\spl \dislag, i\partial_{\bflow} \dislag' \spr.
\end{align}

\subsection{The variational formulation}\label{subsec:varform}

If we test \eqref{eq:PDE} with smooth test functions $(\dislag',\graveu')$ and exploit \eqref{eq:decay}, \eqref{eq:cont-equi} and \eqref{eq:IntbyPartsBflow}, we obtain
\begin{align}\label{eq:variationaleq}
\ses\big((\dislag,\graveu),(\dislag',\graveu')\big) = \spl \boldf,\dislag' \spr
\end{align}
with the sesquilinear form
\begin{align}\label{eq:sesq-a}
\begin{split}
\ses\big((\dislag,\graveu),(\dislag',\graveu')\big)
&:=\spl \sounds^2 \div \dislag, \div\dislag' \spr
+\spl \dens^{-1} \nabla\pres\cdot \dislag, \div\dislag' \spr
+\spl \div\dislag, \dens^{-1} \nabla\pres\cdot \dislag' \spr\\
&-\spl (\omega+i\partial_\bflow+i\angvel\times) \dislag, (\omega+i\partial_\bflow+i\angvel\times) \dislag' \spr\\
&+\spl (\dens^{-1}\hess(\pres)-\hess(\grav))\dislag,\dislag' \spr
-i\omega \spl \damp \dislag, \dislag' \spr\\
&- \spl \nabla\graveu, \dislag' \spr
- \spl \dislag, \nabla\graveu' \spr
+\frac{1}{4\pi G} \spl \nabla\graveu, \nabla\graveu' \spr_{\boldL^2}.
\end{split}
\end{align}
Note that by means of $\qcoeff$ and $\mati$ we can reformulate
\begin{align}\label{eq:sesq-q-mati}
\begin{split}
\ses\big((\dislag,\graveu),(\dislag',\graveu')\big)
&=\spl \sounds^2 (\div \dislag + \qcoeff\cdot\dislag), \div \dislag'  + \qcoeff\cdot\dislag'\spr\\
&-\spl (\omega+i\partial_\bflow+i\angvel\times) \dislag, (\omega+i\partial_\bflow+i\angvel\times) \dislag' \spr\\
&-\spl \mati \dislag,\dislag' \spr-i\omega \spl \damp \dislag, \dislag' \spr\\
&- \spl \nabla\graveu, \dislag' \spr
- \spl \dislag, \nabla\graveu' \spr
+\frac{1}{4\pi G} \spl \nabla\graveu, \nabla\graveu' \spr_{\boldL^2}.
\end{split}
\end{align}
Hence, for Lipschitz domains $D\subset\setR^3$ let
\begin{align*}
\hs(D)&:=\{\dislag \in \boldL_\dens^2(D)\colon \quad\div \dislag \in L_\dens^2(D), \quad \partial_\bflow \dislag \in \boldL_\dens^2(D)\},\\
\spl \dislag,\dislag' \spr_{\hs(D)} &:= \spl \div \dislag, \div \dislag' \spr_D
+\spl \partial_\bflow\dislag, \partial_\bflow\dislag' \spr_D
+\spl \dislag, \dislag' \spr_D,
\end{align*}
and set
\begin{align}\label{eq:hsX}
\hs:=\hs(\setR^3), \qquad \spl \cdot,\cdot \spr_{\hs}:=\spl \cdot,\cdot \spr_{\hs(\setR^3)}.
\end{align}
In addition, for Lipschitz domains $D\subset\setR^3$ with non-trivial boundary let
\begin{align*}
\hs_0(D):=\{\boldu\in\hs \colon \nv\cdot\boldu=0 \text{ on }\partial D\},
\qquad \spl \cdot,\cdot \spr_{\hs_0(D)}:=\spl \cdot,\cdot \spr_{\hs(D)}.
\end{align*}
Under the assumptions of Section~\ref{subsec:assumptions} it follows as in \cite[Lemma~2.1]{HallaHohage:21} that the spaces $\hs$, $\hs(D)$, $\hs_0(D)$ are Hilbert spaces.
The appropriate space for the Eulerian perturbation of gravity $\graveu$ is a bit more technical, because the $L^2(\setR^3)$-norm of $\graveu$ cannot be bounded by means of the sesquilinear form $\ses(\cdot,\cdot)$.
Let
\begin{align}\label{eq:hsgrav}
\begin{split}
\hsgrav&:=\Big\{\graveu\in H^1_\mathrm{loc}(\setR^3) \colon \nabla\graveu\in\boldL^2
\quad\text{and}\quad
\int_{B_{r_1}} \graveu \,\dd\boldx=0\Big\},\\
\spl\graveu,\graveu'\spr_{\hsgrav}&:=\spl\nabla\graveu,\nabla\graveu'\spr_{\boldL^2(\setR^3)}.
\end{split}
\end{align}
By the standard Helmholtz decomposition $\boldG:=\{\boldg\in \boldL^2(\setR^3)\colon \curl\boldg=0\}$ is a closed subspace of $\boldL^2(\setR^3)$ and hence a Hilbert space with the $\boldL^2(\setR^3)$-inner product.
For each $\boldg\in\boldG$ there exists a unique gradient potential $\graveu \in \hsgrav$ such that $\boldg=\nabla\graveu$.
Since the map $\graveu\mapsto\nabla\graveu$ is an isometric isomorphism between $\hsgrav$ and $\boldG$, it follows that $\hsgrav$ is a Hilbert space.
Our specific choice of the Gauge condition, $\int_{B_{r_1}} \graveu \,\dd\boldx=0$, in the definition of $\hsgrav$ will turn out useful in Section~\ref{sec:full}.
Hence, both $\hs$ and $\hsgrav$ are well defined Hilbert spaces and it is straightforward to see that the sesquilinear form $a(\cdot,\cdot)$ is well-defined and bounded on $(\hs\times\hsgrav)\times(\hs\times\hsgrav)$.

\subsection{Common definitions}\label{subsec:commondefs}

We introduce some common functional framework.
For some notions it is more convenient to work with operators instead of sesquilinear forms.
Thus, for generic Hilbert spaces $(\hsgeneric, \spl\cdot,\cdot\spr_{\hsgeneric})$, $(\hsgeneric_1, \spl\cdot,\cdot\spr_{\hsgeneric_1})$, $(\hsgeneric_2, \spl\cdot,\cdot\spr_{\hsgeneric_2})$ we introduce the space $\blo(\hsgeneric_1,\hsgeneric_2)$ of bounded linear operators from $\hsgeneric_1$ to $\hsgeneric_2$ and set $\blo(\hsgeneric):=\blo(\hsgeneric,\hsgeneric)$.
For $\tilde\op\in\blo(\hsgeneric_1,\hsgeneric_2)$ we call $\tilde\op^*\in\blo(\hsgeneric_2,\hsgeneric_1)$ its adjoint, which is defined through $\spl \ugeneric,\tilde\op^*\ugeneric' \spr_{\hsgeneric_1} = \spl \tilde\op\ugeneric,\ugeneric' \spr_{\hsgeneric_2}$ for all $\ugeneric\in\hsgeneric_1, \ugeneric'\in\hsgeneric_2$.
We denote sesquilinear forms with lower case letters and operators with upper case letters.
For a bounded sesquilinear form $\tilde \ses(\cdot,\cdot)$ let $\tilde \op\in \blo(\hsgeneric)$ be its Riesz representation, which is characterized by the relation
\begin{align}\label{eq:defRieszRep}
\spl \tilde \op \ugeneric, \ugeneric' \spr_\hsgeneric = \tilde \ses(\ugeneric,\ugeneric') \quad\text{for all}\quad \ugeneric,\ugeneric'\in\hsgeneric.
\end{align}
Vice-versa for $\tilde\op\in\blo(\hsgeneric)$ let $\tilde\ses(\cdot,\cdot)$ be the bounded sesquilinear form defined by the left-hand-side of \eqref{eq:defRieszRep}.
The tildes in the previous definition were merely used to prevent a confusion with the sesquilinear form $\ses(\cdot,\cdot)$ defined in \eqref{eq:sesq-a}.
The variational equation \eqref{eq:variationaleq} can now be reformulated as operator equation
\begin{align*}
A(\dislag,\graveu)=(\boldf,0).
\end{align*}
(with $\boldf\in\hs$).

\begin{definition}
We say that $\tilde\op\in\blo(\hsgeneric)$ is coercive, if $\inf_{\ugeneric\in\hsgeneric\setminus\{0\}} |\spl \tilde\op \ugeneric,\ugeneric \spr_\hsgeneric|/\|\ugeneric\|^2_\hsgeneric$ $>0$.
We say that $\tilde\op\in\blo(\hsgeneric)$ is weakly coercive, if there exists compact $K\in\blo(\hsgeneric)$ such that $\tilde\op+K$ is coercive.
We say that $\tilde\op\in\blo(\hsgeneric)$ is (weakly) $T$-coercive, if $T\in\blo(\hsgeneric)$ is bijective and $T^*\tilde\op$ is (weakly) coercive.
The same coercivity properties are also attributed to the associated sesquilinear form $\tilde\ses$ defined by \eqref{eq:defRieszRep}.
\end{definition}
The following proposition follows easily from the Lax-Milgram lemma and Riesz theory:
\begin{proposition}\label{prop:Fredholm}
If $\tilde\op$ is weakly $T$-coercive, then $\tilde\op$ is a Fredholm operator with index zero. 
\end{proposition}

Recall that a vector space $\hsgeneric$ is called the direct algebraic sum of subspaces $\hsgeneric_1,\dots,\hsgeneric_N\subset \hsgeneric$, denoted by 
\begin{align}\label{eq:alg_sum}
\hsgeneric=\bigoplus_{n=1,\dots,N}\hsgeneric_n
\end{align}
if each element $y\in \hsgeneric$ has a unique representation of the form $y=\sum_{n=1}^Ny_n$ with $y_n\in 
\hsgeneric_n$.
We refer to \eqref{eq:alg_sum} as algebraic decomposition of $\hsgeneric$. 
Note that there exist associated projection operators $P_{\hsgeneric_n}:\hsgeneric\to \hsgeneric_n$, $y\mapsto y_n$ with 
$\ran P_{\hsgeneric_n}=\hsgeneric_n$ and $\ker P_{\hsgeneric_n} = \bigoplus_{m=1,\dots,N,m\neq n}\hsgeneric_m$.  
\begin{definition}\label{def:topdecomp}
An algebraic decomposition \eqref{eq:alg_sum} of a Hilbert space $\hsgeneric$ is called a \emph{topological decomposition}, denoted by $\bigoplus^{\calT}$, if all associated projection  operators $P_{\hsgeneric_n}$ are continuous. 
\end{definition}
Note that in a topological decomposition all subspaces $\hsgeneric_n = \bigcap_{m\neq n} \ker P_{\hsgeneric_m}$ are closed and the norms $\|y\|_\hsgeneric$ and $\sqrt{\sum_{n=1}^N \|P_{\hsgeneric_n}y\|^2_\hsgeneric}$ are equivalent.

\section{Cowling approximation}\label{sec:cowling}

A common approximation to~\eqref{eq:PDE}, the so-called Cowling approximation, is to set $\graveu\equiv0$ in \eqref{eq:PDE-dislag} and to discard Equation \eqref{eq:PDE-gravlag} together with the decay condition for $\graveu$, i.e.\
\begin{align*}
\begin{split}
-\nabla \big(\sounds^2 \dens\div \dislag
+\nabla \pres\cdot\dislag\big)
+\nabla \pres \div \dislag
+\hess(\pres)\dislag
& \phantom{= \dens\boldf \quad\text{in }\setR^3,}\\
-\dens\hess(\grav)\dislag
- \dens\big(\omega +i\partial_{\bflow}+i\angvel\times \big)^2 \dislag
-i\omega\damp\dens\dislag
&= \dens\boldf \quad\text{in }\setR^3,
\end{split}\\
\dislag &\in H(\div,\setR^3).
\end{align*}
The corresponding sesquilinearform is
\begin{align}\label{eq:ses-cow}
\sescow(\dislag,\dislag'):=\ses\big((\dislag,0),(\dislag',0)\big).
\end{align}
In this section we analyze $\sescow(\cdot,\cdot)$.
Later on, in Section~\ref{sec:full} we will generalize our results to $\ses(\cdot,\cdot)$.
First, we note the injectivity of $\opAcow$.
\begin{lemma}\label{lem:inj-cow}
Let the assumptions of Section~\ref{subsec:assumptions} be satisfied and $\omega\neq0$.
Then $\opAcow$ is injective.
\end{lemma}
\begin{proof}
Let $\dislag\in\ker\opAcow$. Then 
\begin{align*}
0=\big|\Im \big(\sescow(\dislag,\dislag)\big)\big|
= |\omega| \spl \damp \dislag,\dislag \spr
\geq |\omega|\dampl \|\dislag\|^2_{\boldL^2_\dens}
\end{align*}
and hence $\dislag=0$.
\end{proof}
Our goal is to prove the bijectivity of $\opAcow$.
To succeed, it remains to show that $\opAcow$ is weakly T-coercive.
This result was achieved in \cite{HallaHohage:21} for bounded domains $D$ together with a boundary condition on $\partial D$.
Therein a major ingredient was a topological decomposition $\hs_0(D)=\hsV\oplus^\calT\hsW\oplus^\calT\hsZ$ such that $\hsV\subset\boldH^1(D)$, and hence with the compact embedding $\hsV\hookrightarrow\boldL^2(D)$.
However, for the unbounded domain $\setR^3$ the embedding $\boldH^1(\setR^3)\hookrightarrow\boldL^2(\setR^3)$ is no longer compact, and thus an analogous topological decomposition of $\hs(\setR^3)$ on a global level does not suffice to reproduce the former results.
On the other hand, in the exterior part $B_{r_1}^c$ it holds $\bflow\equiv0$ and hence for $\dislag\in\hs$ with $\supp\dislag\subset B_{r_1}^c$ it follows with \eqref{eq:ses-cow} and \eqref{eq:sesq-q-mati} that
$\sescow(\dislag,\dislag)=\|\sounds(\div+\qcoeff\cdot)\dislag\|_{L^2_\dens}^2-\spl (\matii+i\omega\damp)\dislag,\dislag\spr$.
Since the numerical range of the matrix $\matii(\boldx)+i\omega\damp(\boldx) I_{3\times3}$ is contained in a fixed closed salient sector in the upper half plane for all $\boldx\in B_{r_1}^c$, it follows that $\sescow(\cdot,\cdot)$ is coercive in the exterior domain $B_{r_1}^c$.
Thus, our approach is to combine a topological decomposition based in the interior together with the coercivity in the exterior to a unified analysis.
However, it is a delicate matter to marry these two separate ideas.
As preparation we introduce in the next theorem a topological decomposition of $\hs$ which is similar, but different to the one in \cite[Theorem~3.5]{HallaHohage:21}.

\begin{theorem}\label{thm:decomposition}
Let $\dens$, $\bflow$, $\qcoeff$ and $r_1$ be as in Section~\ref{subsec:assumptions}.
Let $r_2>r_1$.
Then $\hs$ admits a topological decomposition  
\begin{align*}
\hs = \hsV \oplus^\calT \hsW \oplus^\calT \hsZ
\end{align*}
with the following properties:
\begin{enumerate}
\item\label{it:Vcompactembedding}
$\hsV\subset\{\dislag|_{B_{r_2}^c}=0, \dislag|_{B_{r_2}}=\nabla\compvp\colon
\compvp\in H^2(B_{r_2})\mbox{ with }\frac{\partial \compvp}{\partial \nv}=0 \mbox{ on }\partial B_{r_2}\}$ is compactly embedded in $\boldL^2$.
\item\label{it:Wkernel}
$\hsW = \{\dislag \in \hs \colon \div \dislag + \qcoeff\cdot \dislag = 0 \text{ in } B_{r_2} \} $.
\item\label{it:Zfinite}
$\hsZ$ is finite-dimensional. 
\end{enumerate}
Moreover, for each $\zeta\in W^{1,\infty}(B_{r_2})$ there exists a compact operator $K_\zeta\in\blo(\hs)$ such that
\begin{align}\label{eq:multiplier}
\|\zeta\div\compv\|_{L^2(B_{r_2})}^2=
\|\zeta\nabla\compv\|_{(L^2(B_{r_2}))^{3x3}}^2
+\spl K_\zeta\compv,\compv \spr_{\hs}
\end{align}
for all $\compv\in\hsV$.
\end{theorem}
\begin{proof}
We apply the same construction as in \cite[Theorem~3.5]{HallaHohage:21}.
Hence we only sketch the main ideas and focus on the minor points which diverge from \cite[Theorem~3.5]{HallaHohage:21}.
For given $\dislag\in\hs$ we consider the problem to find $\compvp\in H^1(B_{r_2})$ such that
\begin{subequations}\label{eq:EqForv0}
\begin{align}
\Delta\compvp+\qcoeff\cdot\nabla\compvp&=\div\dislag+\qcoeff\cdot\dislag\phantom{0}\hspace{-2mm} \quad\text{in }B_{r_2},\\
\nv\cdot\nabla\compvp &=0\phantom{\div\dislag+\qcoeff\cdot\dislag}\hspace{-2mm} \quad\text{on }\partial B_{r_2}.
\end{align}
\end{subequations}
Since $\qcoeff\cdot\nabla$ is a (low order) perturbation of $\Delta$, it is not guaranteed that \eqref{eq:EqForv0} admits a unique solution.
However, the perturbation is compact and hence there exists a finite dimensional subspace $\hsZ\subset\hs$ with a projection $\PZ\in L(\hs,\hsZ)$ onto $\hsZ$ such that
\begin{subequations}\label{eq:EqForv0Z}
\begin{align}
\Delta\compvp+\qcoeff\cdot\nabla\compvp&=\div\big((1-\PZ)\dislag\big)+\qcoeff\cdot(1-\PZ)\dislag\phantom{0}\hspace{-2mm} \quad\text{in }B_{r_2},\\
\nv\cdot\nabla\compvp &=0\phantom{\div\big((1-\PZ)\dislag\big)+\qcoeff\cdot(1-\PZ)\dislag}\hspace{-2mm} \quad\text{on }\partial B_{r_2}.
\end{align}
\end{subequations}
admits a unique solution $\compvp\in H^1(B_{r_2})$ with $\|\compvp\|_{H^1(B_{r_2})}\lesssim \|\dislag\|_\hs$.
Due to $\dens\in L^\infty$ it holds $\nabla\compvp\in\boldL^2_\dens(B_{r_2})$.
Since $\Delta\compvp\in L^2(B_{r_2})$, $\nv\cdot\nabla\compvp=0$ on $\partial B_{r_2}$, and because $B_{r_2}$ is a $C^{1,1}$-domain it follows with convenient regularity theory (e.g.\ \cite[Theorems 2.3.3.2, 2.4.2.7]{Grisvard:85}) that $\compvp\in H^2(B_{r_2})$ and
\begin{align*}
\|\compvp\|_{H^2(B_{r_2})}\lesssim \|\compvp\|_{H^1(B_{r_2})}+\|\Delta\compvp\|_{L^2(B_{r_2})} \lesssim \|\dislag\|_\hs.
\end{align*}
Thus $\partial_\bflow\nabla\compvp\in \boldL^2_\dens(B_{r_2})$ and hence $\nabla\compvp\in \hs(B_{r_2})$ with $\|\nabla\compvp\|_{\hs(B_{r_2})} \lesssim \|\dislag\|_\hs$.
Let $\compv$ be the continuation by zero of $\nabla\compvp$ to $B_{r_2}^c$.
Then $\compv\in\hs(B_{r_2}^c)$.
Since the normal traces of $\compv|_{B_{r_2}}$ and $\compv|_{B_{r_2}^c}$ both vanish at $\partial B_{r_2}$ it follows $\div\compv\in L^2_\dens$.
Since $\bflow$ vanishes in a neighborhood of $\partial B_{r_2}$ it also follows $\partial_\bflow\compv\in\boldL^2_\dens$.
Thus $\compv\in\hs$ and $\|\compv\|_\hs\lesssim\|\dislag\|_\hs$.
Hence for $\PV\dislag:=\compv$ it holds $\PV\in L(\hs)$ and we set $\hsV:=\ran\PV$.
Finally, we set $\PW\dislag:=\dislag-(\PV+\PZ)\dislag$ and it follows $(\div+\qcoeff\cdot)\PW\dislag=0$ in $B_{r_2}$ from \eqref{eq:EqForv0Z}.
Formula \eqref{eq:multiplier} follows the same way as in \cite[Theorem~3.5]{HallaHohage:21}.
\end{proof}
To construct a suitable $T$-operator consider $\mu\in C^\infty(\setR)$ with
\begin{enumerate}
 \item $\mu(r)=0$ for $r\leq r_1$,
 \item $\mu$ is non-decreasing,
 \item $\mu(r)=\mu_*$ for $r\geq r_2$ with constant $\mu_*\in(0,\pi)$.
\end{enumerate}
For such $\mu$ let $\sigma(x):=e^{i\mu(|x|)\sign\omega}$ and $\sigma_*:=e^{i\mu_*\sign\omega}$.
Then a natural candidate for $T$ is $T=\PV-\sigma\PW+\PZ$.
Although it holds $\sigma\PW\dislag\notin\hsW$, which is undesirable for our analysis. 
Hence we introduce a slight modification.
For given $\compw\in\hsW$ consider the problem to find $\hat v_0\in H^1(A_{r_1,r_2})$ such that
\begin{align*}
(\div+\qcoeff\cdot)(\nabla \hat v_0-\sigma\compw)&=0 \quad\text{in } A_{r_1,r_2},\\
\nv\cdot\nabla \hat v_0&=0 \quad\text{on }\partial A_{r_1,r_2}.
\end{align*}
Since this equation is weakly coercive, we can find a finite dimensional subspace $\hsW_0\subset\hsW$ and a projection $K_1\in L(\hsW,\hsW_0)$ onto $\hsW_0$ such that the problem to find $\hat v_0\in H^1(A_{r_1,r_2})$ with
\begin{align*}
(\div+\qcoeff\cdot)(\nabla \hat v_0-\sigma\compw+\sigma K_1\compw)&=0 \quad\text{in } A_{r_1,r_2},\\
\nv\cdot\nabla \hat v_0&=0 \quad\text{on }\partial A_{r_1,r_2},
\end{align*}
admits a unique solution $\hat v_0$.
With convenient regularity theory (e.g.\ \cite[Theorems 2.3.3.2, 2.4.2.7]{Grisvard:85}) it follows that $\hat v_0 \in H^2(A_{r_1,r_2})$ and $\|\hat v_0\|_{H^2(A_{r_1,r_2})} \lesssim \|\boldu\|_\hs$.
Let $\hat \boldv$ be the continuation of $\nabla \hat v_0$ to $\setR^3\setminus A_{r_1,r_2}$ by zero.
It follows that $\hat \boldv\in \hs$.
Hence let
\begin{align*}
T_\hsW' \boldw:=\sigma\boldw-\hat\boldv-\sigma K_1\boldw.
\end{align*}
It follows from the definition of $\hat\boldv$ that $(\div+\qcoeff\cdot)T_\hsW'\boldw=0$ in $A_{r_1,r_2}$.
Since $\sigma=1$ in $B_{r_1}$, $K_1$ maps into $\hsW_0\subset\hsW$ and $\hat\boldv=0$ in $B_{r_1}$ it follows $(\div+\qcoeff\cdot)T_\hsW'\boldw=0$ in $B_{r_1}$ too.
Hence $T_\hsW'\in L(\hsW)$.
However, is $T_\hsW'$ invertible?
To shed light on this question we compute
\begin{align*}
\spl T_\hsW'\boldw,\boldw \spr_\hs&=
\|\boldw\|_{\hs(B_{r_1})}^2
-\spl K_1\boldw,\boldw \spr_{\hs(B_{r_1})}\\
&+\spl \sigma\boldw,\boldw \spr_{A_{r_1,r_2}}
-\spl \hat\boldv,\boldw \spr_{A_{r_1,r_2}}
-\spl \sigma K_1\boldw,\boldw \spr_{A_{r_1,r_2}}\\
&+\spl \sigma\qcoeff\cdot\boldw,\qcoeff\cdot\boldw \spr_{A_{r_1,r_2}}
-\spl \qcoeff\cdot\hat\boldv,\qcoeff\cdot\boldw \spr_{A_{r_1,r_2}}
-\spl \qcoeff\cdot \sigma K_1\boldw,\qcoeff\cdot\boldw \spr_{A_{r_1,r_2}}\\
&+\sigma_*\|\boldw\|_{\hs(B_{r_2}^c)}^2
-\sigma_* \spl K_1\boldw,\boldw \spr_{\hs(B_{r_2}^c)}.
\end{align*}
Thereat we exploited that $\sigma$ is constant in $\setR^3\setminus A_{r_1,r_2}$, $\hat\boldv=0$ in $\setR^3\setminus A_{r_1,r_2}$, $\bflow=0$ in $B_{r_1}^c$ and $\div\boldw=-\qcoeff\cdot\boldw$ in $B_{r_2}$ for $\boldw\in\hsW$.
All terms which involve $K_1$ are compact, because $K_1$ is compact.
Also, all terms which involve $\hat\boldv$ are compact due to the compact embedding $\boldH^1(A_{r_1,r_2})\hookrightarrow \boldL^2(A_{r_1,r_2})$.
Again with  $\div\boldw=-\qcoeff\cdot\boldw$ in $B_{r_2}$ for $\boldw\in\hsW$, we obtain that the remainder
\begin{align*}
\|\boldw\|_{\hs(B_{r_1})}^2
+\spl \sigma\boldw,\boldw \spr_{A_{r_1,r_2}}
+\spl \sigma\qcoeff\cdot\boldw,\qcoeff\cdot\boldw \spr_{A_{r_1,r_2}}
+\sigma_*\|\boldw\|_{\hs(B_{r_2}^c)}^2
\end{align*}
is coercive due to $\sigma(\boldx)=e^{i\mu(|\boldx|)\sign\omega}$, $\mu\in[0,\mu_*]$ and $\mu_*\in(0,\pi)$.
Thus $T_\hsW'$ is weakly coercive.
Hence there exists a compact operator $K_2\in L(\hsW)$ such that
\begin{align*}
T_\hsW:=T_\hsW'+K_2
\end{align*}
is coercive and hence bijective.
Now we can define
\begin{align*}
T:=\PV-T_\hsW\PW+\PZ,
\end{align*}
which is bijective with inverse $T^{-1}=\PV-T_\hsW^{-1}\PW+\PZ$.
In preparation of Theorem~\ref{thm:wTc-cow} we introduce a technical quantity $\theta$.
To this end let the function $\arg$ take values in $(-\pi,\pi]$ and recall that the numerical range of a matrix $M\in\setC^{3\times3}$ is defined by
\begin{align*}
\numran M := \{ \xi^H M \xi\colon \xi\in\setC^3, |\xi|=1\},
\end{align*}
with the Euclidean vector norm $|\cdot|$.
Then let
\begin{align}\label{eq:theta}
\theta&:=\max\Big\{0, \sup_{\boldx\in B_{r_2}} |\sup \arg \numran \big(i\omega\damp(\boldx) I_{3\times3}+\mati(\boldx)\big)|-\pi/2\Big\}
\end{align}
Note that the numerical range of a matrix is a set and the inner supremum/infimum in definitions like for $\theta$ is over this set.
Further, we can estimate
\begin{align*}
\theta \leq \arg\big( i|\omega|\dampl-\|M_{\mati}\|_{L(\boldL^2_\dens(B_{r_2}))} \big) -\pi/2< \pi/2.
\end{align*}

\begin{theorem}\label{thm:wTc-cow}
Let the assumptions of Section~\ref{subsec:assumptions} be satisfied.
Let $r_2>r_1$ be such that $\sounds, \dens \in W^{1,\infty}(B_{r_2})$ and $\omega\neq0$.
If
\begin{align*}
\|\sounds^{-1}\bflow\|_{\boldL^\infty}^2<\frac{1}{1+\tan^2\theta},
\end{align*}
then there exists $\mu$ such that $A_\mathrm{Cow}$ is weakly $T$-coercive.
\end{theorem}
\begin{proof}
For $\dislag,\dislag'\in\hs$ we use the notation
\begin{align*}
&\compv:=\PV\dislag, &&\compw:=\PW\dislag, && \compz:=\PZ\dislag,\\
&\compv':=\PV\dislag',&&\compw':=\PW\dislag', &&\compz':=\PZ\dislag'
\end{align*}
such that $\dislag = \compv+\compw+\compz$ and $\dislag'= \compv'+\compw'+\compz'$. 
We introduce a sufficiently small parameter $\tau\in(0,\pi/2-\theta)$ which will be specified later on.
Subsequently we set $\beta:=\mu_*-\pi/2+\theta+\tau$.
We choose $\mu$ such that $\beta\in(0,\pi/2)$ and
\begin{align}\label{eq:EstBeta}
\begin{split}
0&<\Re \big( -ie^{-i\sign\omega\beta} (i|\omega|\dampl+\sign\omega\|M_{\matii}\|_{L(\boldL^2_\dens)} \big)\\
&=\Re \big( e^{-i\beta} (|\omega|\dampl-i\|M_{\matii}\|_{L(\boldL^2_\dens)} \big).
\end{split}
\end{align}
Note that such a choice is possible, and $\beta\in(0,\pi/2)$ implies $\mu_*<\pi$ (which is required for the construction of $T_\boldW$).

\emph{definition of $A_1$ and $A_2$:}
We split $T^*\opAcow=A_1+A_2$ into a coercive operator $A_1$ and a compact operator $A_2$.
To this end we introduce an additional small parameter $\delta>0$, which will be specified later on.
Let $K_{\sounds\sqrt{\dens}}$ be as in Theorem~\ref{thm:decomposition}.
Then we define $A_1$ and $A_2$ by
\begin{align*}
\spl A_1^\mathrm{int}\dislag,\dislag'\spr_\hs &:=
\spl \sounds^2\div \compv,\div \compv'\spr
-\spl i\partial_\bflow \compv,i\partial_\bflow \compv' \spr
+\spl \compv,\compv'\spr
+\frac{1}{4\delta}\spl K_{\sounds\sqrt{\dens}}\compv,K_{\sounds\sqrt{\dens}}\compv' \spr_\hs\\
&-\spl (\omega+i\partial_\bflow+i\angvel\times)\compw,i\partial_\bflow \compv' \spr
+\spl i\partial_\bflow \compv,(\omega+i\partial_\bflow+i\angvel\times)\compw' \spr
\\
&+\spl (\omega+i\partial_\bflow+i\angvel\times)\compw,(\omega+i\partial_\bflow+i\angvel\times)\compw' \spr_{B_{r_1}},\\
&+\spl (i\omega\damp+\mati)\compw,\compw'\spr_{B_{r_1}}\\
\spl A_1\dislag,\dislag'\spr_\hs &:=
\spl A_1^\mathrm{int}\dislag,\dislag'\spr_\hs\\
&+\spl \ol{\sigma}(i\omega\damp+\matii)\compw,\compw'\spr_{B_{r_1}^c}\\
&-\ol{\sigma_*} \spl \sounds^2(\div+\qcoeff\cdot)\compw,(\div+\qcoeff\cdot)\compw'\spr_{B_{r_2}^c}\\
&-\ol{\sigma_*} \spl \compz,\compz' \spr
\end{align*}
and
\begin{align*}
\spl \op_2\dislag,\dislag'\spr_\hs &=
-\spl \compv,\compv'\spr
-\frac{1}{4\delta}\spl K_{\sounds\sqrt{\dens}}\compv,K_{\sounds\sqrt{\dens}}\compv' \spr_\hs
+\ol{\sigma_*}\spl \compz,\compz' \spr_\hs\\
&+\spl \div\compv,\nabla\pres\cdot\compv' \spr_{L^2}
+\spl \nabla\pres\cdot\compv,\div\compv' \spr_{L^2}
+\spl \sounds^{-2}\dens^{-1}\; \nabla\pres\cdot\compv, \nabla\pres\cdot\compv' \spr_{L^2}\\
&-\spl i\partial_\bflow \compv,(\omega+i\angvel\times)\compv' \spr
-\spl (\omega+i\angvel\times)\compv,i\partial_\bflow \compv' \spr
-\spl (i\omega\damp+\matii) \compv,\compv' \spr\\
&-\spl (\omega+i\partial_\bflow+i\angvel\times)\compw,(\omega+i\angvel\times) \compv' \spr
+\spl (\omega+i\angvel\times) \compv,(\omega+i\partial_\bflow+i\angvel\times)T'_\boldW\compw' \spr\\
&-\spl (i\omega\damp+\mati) \compw,\compv' \spr
+\spl (i\omega\damp+\mati) \compv,T'_\boldW\compw' \spr\\
&-\spl i\partial_\bflow \compv,(\omega+i\partial_\bflow+i\angvel\times)\sigma K_1\compw' \spr\\
&-\spl (\omega+i\partial_\bflow+i\angvel\times)\compw,(\omega+i\partial_\bflow+i\angvel\times)K_1\compw' \spr_{B_{r_1}}
-\spl (i\omega\damp+\mati)\compw,K_1\compw'\spr_{B_{r_1}}\\
&-\spl \ol{\sigma}(i\omega\damp+\matii)\compw,K_1\compw'\spr_{B_{r_1}^c}
+\ol{\sigma_*} \spl \sounds^2(\div+\qcoeff\cdot)\compw,(\div+\qcoeff\cdot)K_1\compw'\spr_{B_{r_2}^c}\\
&-\spl (i\omega\damp+\matii)\compw,\hat\boldv'\spr_{A_{r_1,r_2}}\\
&-\sescow(\boldv+\boldw,K_2\boldw')\\
&+\sescow(\compv+\compw,\compz')
+\sescow(\compz,\compv'-T_\boldW\compw')
+\sescow(\compz,\compz')
\end{align*}
for all $\dislag,\dislag'\in\hs$.
To see that indeed $T^*\opAcow=A_1+A_2$ consider the following.
The terms
\begin{align*}
\spl \compv,\compv'\spr
+\frac{1}{4\delta}\spl K_{\sounds\sqrt{\dens}}\compv,K_{\sounds\sqrt{\dens}}\compv' \spr_\hs
-\ol{\sigma_*} \spl \compz,\compz' \spr
\end{align*}
are compact and added into $A_1$ to guarantee the coercivity of $A_1$.
Subsequently these terms are added with a reverse sign into $A_2$ to sustain the identity $T^*\opAcow=A_1+A_2$.
All remaining terms which involve $\boldz$ or $\boldz'$
\begin{align*}
\sescow(\compv+\compw,\compz')+\sescow(\compz,\compv'-T_\boldW\compw')
+\sescow(\compz,\compz')
\end{align*}
are put into $A_2$.
With the definitions of $\qcoeff$ and $\matii$ in \eqref{eq:q-mati} and \eqref{eq:matii} respectively a complete expansion of $\sescow(\boldv,\boldv')$ yields
\begin{align*}
\sescow(\boldv,\boldv')&=
\spl \sounds^2\div \compv,\div \compv'\spr
-\spl i\partial_\bflow \compv,i\partial_\bflow \compv' \spr
+\spl \div\compv,\nabla\pres\cdot\compv' \spr_{L^2}\\
&+\spl \nabla\pres\cdot\compv,\div\compv' \spr_{L^2}+\spl \sounds^{-2}\dens^{-1}\; \nabla\pres\cdot\compv, \nabla\pres\cdot\compv' \spr_{L^2}\\
&-\spl i\partial_\bflow \compv,(\omega+i\angvel\times)\compv' \spr
-\spl (\omega+i\angvel\times)\compv,i\partial_\bflow \compv' \spr
-\spl (i\omega\damp+\matii) \compv,\compv' \spr,
\end{align*}
and these terms are split into $A_1$ and $A_2$.
Next we apply $T_\boldW=T_\boldW'+K_2$ and all remaining terms involving $K_2$
\begin{align*}
-\sescow(\boldv+\boldw,K_2\boldw)
\end{align*}
are put into $A_2$.
Further, we use representation \eqref{eq:sesq-q-mati} and since $\boldw, T_\boldW'\boldw\in\hsW$ the term
\begin{align*}
\spl \sounds^2 (\div+\qcoeff\cdot)\boldw,\boldv \spr_{B_{r_2}}
-\spl \sounds^2 \boldv, (\div+\qcoeff\cdot)T'_\hsW\boldw' \spr_{B_{r_2}}
\end{align*}
vanishes.
Since $\boldv=0$ in $B_{r_2}^c$ the remaining terms of
$\sescow(\boldw,\boldv')+\sescow(\boldv,-T'_\boldW\boldw')$ read
\begin{align}\label{eq:mixed-terms}
\begin{split}
&-\spl (\omega+i\partial_\bflow+i\angvel\times)\compw,i\partial_\bflow \compv' \spr
+\spl i\partial_\bflow \compv,(\omega+i\partial_\bflow+i\angvel\times)T'_\boldW\compw' \spr\\
&-\spl (\omega+i\partial_\bflow+i\angvel\times)\compw,(\omega+i\angvel\times) \compv' \spr
+\spl (\omega+i\angvel\times) \compv,(\omega+i\partial_\bflow+i\angvel\times)T'_\boldW\compw' \spr\\
&-\spl (i\omega\damp+\mati) \compw,\compv' \spr
+\spl (i\omega\damp+\mati) \compv,T'_\boldW\compw' \spr.
\end{split}
\end{align}
We use $T_\boldW'\boldw'=\sigma\boldw'-\hat\boldv'-\sigma K_1\boldw'$, $\bflow=0$ in $B_{r_1}^c$, $\hat\boldv'=0$ and $\sigma=1$ in $B_{r_1}$ to obtain
\begin{align*}
\spl i\partial_\bflow \compv,(\omega+i\partial_\bflow+i\angvel\times)T'_\boldW\compw' \spr
=\spl i\partial_\bflow \compv,(\omega+i\partial_\bflow+i\angvel\times)\compw' \spr
-\spl i\partial_\bflow \compv,(\omega+i\partial_\bflow+i\angvel\times)\sigma K_1\compw' \spr.
\end{align*}
The part
\begin{align*}
&-\spl (\omega+i\partial_\bflow+i\angvel\times)\compw,i\partial_\bflow \compv' \spr
+\spl i\partial_\bflow \compv,(\omega+i\partial_\bflow+i\angvel\times)\compw' \spr
\end{align*}
of \eqref{eq:mixed-terms} is put into $A_1$ and the remainder of \eqref{eq:mixed-terms} is put into $A_2$.
It remains to discuss $\sescow(\boldw,-T'_\boldW\boldw')$.
The term $-\spl \sounds^2(\div+\qcoeff\cdot)\boldw,(\div+\qcoeff\cdot)T'_\boldW\boldw \spr_{B_{r_2}}$ vanishes due to $\boldw\in\hsW$.
With \eqref{eq:sesq-q-mati} and $\bflow=0$ in $B_{r_1}^c$ we obtain
\begin{align*}
\sescow(\boldw,-T'_\boldW\boldw')&=
\spl (\omega+i\partial_\bflow+i\angvel\times)\compw,(\omega+i\partial_\bflow+i\angvel\times)T'_\boldW\compw' \spr_{B_{r_1}},\\
&+\spl (i\omega\damp+\mati)\compw,T'_\boldW\compw'\spr_{B_{r_1}}\\
&+\spl (i\omega\damp+\matii)\compw,T'_\boldW\compw'\spr_{B_{r_1}^c}\\
&-\spl \sounds^2(\div+\qcoeff\cdot)\compw,(\div+\qcoeff\cdot)T'_\boldW\compw'\spr_{B_{r_2}^c}
\end{align*}
With $T_\boldW'\boldw'=\sigma\boldw'-\hat\boldv'-\sigma K_1\boldw'$, $\hat\boldv'=0$ in $A_{r_1,r_2}^c$, $\sigma=1$ in $B_{r_1}$ and $\sigma=\sigma_*$ in $B_{r_2}^c$ we split the former into a part incorporated into $A_1$
\begin{align*}
&\spl (\omega+i\partial_\bflow+i\angvel\times)\compw,(\omega+i\partial_\bflow+i\angvel\times)\compw' \spr_{B_{r_1}}
+\spl (i\omega\damp+\mati)\compw,\compw'\spr_{B_{r_1}}\\
&+\spl \ol{\sigma}(i\omega\damp+\matii)\compw,\compw'\spr_{B_{r_1}^c}
-\ol{\sigma_*} \spl \sounds^2(\div+\qcoeff\cdot)\compw,(\div+\qcoeff\cdot)\compw'\spr_{B_{r_2}^c}
\end{align*}
and a part incorporated into $A_2$
\begin{align*}
&-\spl (\omega+i\partial_\bflow+i\angvel\times)\compw,(\omega+i\partial_\bflow+i\angvel\times)K_1\compw' \spr_{B_{r_1}}
-\spl (i\omega\damp+\mati)\compw,K_1\compw'\spr_{B_{r_1}}\\
&-\spl \ol{\sigma}(i\omega\damp+\matii)\compw,K_1\compw'\spr_{B_{r_1}^c}
+\ol{\sigma_*} \spl \sounds^2(\div+\qcoeff\cdot)\compw,(\div+\qcoeff\cdot)K_1\compw'\spr_{B_{r_2}^c}\\
&-\spl (i\omega\damp+\matii)\compw,\hat\boldv'\spr_{A_{r_1,r_2}}.
\end{align*}
Thus indeed $T^*\opAcow=A_1+A_2$.

\emph{compactness of $A_2$:}
The operator $A_2$ is compact, because of the compact embedding $\hsV\hookrightarrow \boldL^2$ (see Theorem~\ref{thm:decomposition}, Item~\ref{it:Vcompactembedding}), because of the compactness of $K_{\sounds\sqrt{\dens}}$ (see Theorem~\ref{thm:decomposition}), due to the finite dimension of $\hsZ$ (see Theorem~\ref{thm:decomposition}, Item~\ref{it:Zfinite}), because of the compact embedding $\boldH^1(A_{r_1,r_2})\hookrightarrow \boldL^2(A_{r_1,r_2})$ for $\hat\boldv'$ and because $K_1$, $K_2$ are compact.

\emph{coercivity of $A_1$ (1st part):}
To prove that $A_1$ is coercive we estimate
\begin{align*}
\frac{1}{\cos(\theta+\tau)}
\Re \left( e^{-i(\theta+\tau)\sign \omega} \spl A_1\dislag,\dislag \spr_\hs \right).
\end{align*}
To treat the $A_1^\mathrm{int}$ we proceed as in \cite{HallaHohage:21} and compute
\begin{align*}
\frac{1}{\cos(\theta+\tau)}
\Re\Big( e^{-i(\theta+\tau)\sign\omega} &\spl \op_1^\mathrm{int} \dislag,\dislag \spr_\hs \Big) =\\
&\|\sounds\div \compv\|^2_{L^2_\dens}
-\|\partial_\bflow \compv \|^2_{\boldL^2_\dens}
+\|\compv\|_{\boldL^2_\dens}
+\frac{1}{4\delta}\| K_{\sounds\sqrt{\dens}}\compv\|_\hs^2\\
&+\|(\omega+i\partial_\bflow+i\angvel\times)\compw \|^2_{\boldL^2_\dens(B_{r_1})}
+\frac{|\omega|\sin\tau}{\cos(\theta+\tau)} \| \sqrt{\damp} \compw \|^2_{\boldL^2_\dens(B_{r_1})}\\
&-2\tan(\theta+\tau)\sign\omega\,\Im\big(\spl i\partial_\bflow \compv,(\omega+i\partial_\bflow+i\angvel)\compw \spr_{\boldL^2_\dens}\big).
\end{align*}
Estimating the last term by the Cauchy-Schwarz inequality and the weighted Young inequality $2ab\leq (1-\epsilon)^{-1}a^2+(1-\epsilon)b^2$ 
with another parameter $\epsilon\in(0,1)$, $a=\tan(\theta+\tau) \|\partial_\bflow \compv\|_{\boldL^2_\dens}$, and 
$b= \|(\omega+i\partial_\bflow+i\angvel\times)\compw\|_{\boldL^2_\dens(B_{r_1})}$ we obtain 
\begin{align}\label{eq:abhierweiter}
\begin{split}
\frac{1}{\cos(\theta+\tau)}
\Re\Big( &e^{-i(\theta+\tau)\sign\omega} \spl \op_1^\mathrm{int} \dislag,\dislag \spr_\hs \Big) \geq\\
&\| \sounds\div \compv\|^2_{L^2_\dens}
-\big(1+(1-\epsilon)^{-1}\tan^2(\theta+\tau)\big)\|\partial_\bflow \compv \|^2_{\boldL^2_\dens}\\
&+\|\compv\|_{\boldL^2_\dens}+\frac{1}{4\delta}\|K_{\sounds\sqrt{\dens}}\compv\|_\hs^2\\
&+\epsilon\| (\omega+i\partial_\bflow+i\angvel\times)\compw \|^2_{\boldL^2_\dens(B_{r_1})}
+\frac{|\omega|\sin\tau}{\cos(\theta+\tau)} \| \sqrt{\damp}\compw \|^2_{\boldL^2_\dens(B_{r_1})}.
\end{split}
\end{align}
Due to Theorem~\ref{thm:decomposition} it holds
\begin{align*}
\|\sounds\sqrt{\dens}\div\compv\|_{L^2}^2
=\|\sounds\sqrt{\dens}\nabla\compv\|_{(L^2)^{3\times 3}}^2
+\spl K_{\sounds\sqrt{\dens}}\compv,\compv\spr_\hs
\end{align*}
for each $\compv\in\hsV$.
Since by assumption the flow is subsonic ($\|\sounds^{-1}\bflow\|_{\boldL^\infty}^2<1/(1+\tan^2\theta)$), we can choose $\epsilon$ and $\tau$ small enough such that
\begin{align*}
0<1-\big(1+(1-\epsilon)^{-1}\tan^2(\theta+\tau)\big)\|\sounds^{-1}\bflow\|_{\boldL^\infty}^2=:C_{\epsilon,\tau,\theta}.
\end{align*}
Bounding $\|\partial_\bflow\compv\|^2_{\boldL^2_\dens}$ by $\|\sounds^{-1}\bflow\|_{\boldL^\infty}^2\|\sounds\sqrt{\dens}\nabla\compv\|_{(L^2)^{3\times3}}^2$
we can estimate
\begin{align}\label{eq:estimateinproof}
\begin{split}
\|\sounds\div\compv\|^2_{L^2_\dens}
&-\big(1+(1-\epsilon)^{-1}\tan^2(\theta+\tau)\big)\|\partial_\bflow\compv\|^2_{\boldL^2_\dens}\\
&\geq \soundsl^2\densl_{r_2} C_{\epsilon,\tau,\theta} |\compv|_{\boldH^1(B_{r_2})}^2 -|\spl K_{\sounds\sqrt{\dens}}\compv,\compv\spr_\hs|\\
&\geq \soundsl^2\densl_{r_2} C_{\epsilon,\tau,\theta} |\compv|_{\boldH^1(B_{r_2})}^2 -\frac{1}{4\delta}\|K_{\sounds\sqrt{\dens}}\compv\|_\hs^2 
-\delta\|\compv\|_\hs^2
\end{split}
\end{align}
with $\soundsl$ as in \eqref{eq:box_constraints} and $\densl_{r_2}$ as in \eqref{eq:densl}.
We continue and estimate by means of \eqref{eq:estimateinproof} the first two lines of the right hand side of \eqref{eq:abhierweiter}
\begin{align*}
\|\sounds\div \compv\|^2_{L^2_\dens}
&-\big(1+(1-\epsilon)^{-1}\tan^2(\theta+\tau)\big)\|\partial_\bflow \compv \|^2_{\boldL^2_\dens}
+\|\compv\|_{\boldL^2_\dens}+\frac{1}{4\delta}\|K_{\sounds\sqrt{\dens}}\compv\|_\hs^2\\
&\geq \soundsl^2\densl C_{\epsilon,\tau,\theta} |\compv|_{\boldH^1(B_{r_2})}^2 +\|\compv\|_{\boldL^2_\dens}^2 -\delta\|\compv\|_\hs^2.
\end{align*}
There exists a constant $C_V>0$ such that
\begin{align*}
\soundsl^2\densl C_{\epsilon,\tau,\theta} |\compv|_{\boldH^1(B_{r_2})}^2 +\|\compv\|_{\boldL^2_\dens}^2 \geq C_V\|\compv\|_\hs^2
\end{align*}
for each $\compv\in\hsV$. Thus
\begin{align}\label{eq:InequProofG1}
\begin{aligned}
\| \sounds\div \compv\|^2_{L^2_\dens}
&-\big(1+(1-\epsilon)^{-1}\tan^2(\theta+\tau)\big)\| \partial_\bflow \compv \|^2_{\boldL^2_\dens}
+\|\compv\|_{\boldL^2_\dens}+\frac{1}{4\delta}\|K_{\sounds\sqrt{\dens}}\compv\|_\hs^2\\
&\geq (C_V-\delta)\|\compv\|_\hs^2.
\end{aligned}
\end{align}
Now we choose $\delta<C_V$.
The third line of the right hand side of \eqref{eq:abhierweiter}
can be estimated using a weighted Young inequality and $\dampl>0$ (with $\dampl$ defined in \eqref{eq:box_constraints}) to obtain
\begin{align}\label{eq:InequProofG2}
\epsilon\|(\omega+i\partial_\bflow+i\angvel\times)\compw\|^2_{\boldL^2_\dens(B_{r_1})}
&+\frac{|\omega|\sin\tau}{\cos(\theta+\tau)}\|\sqrt{\damp}\compw\|^2_{\boldL^2_\dens(B_{r_1})}
\geq C_W(\|\partial_\bflow\compw\|^2_{\boldL^2_\dens(B_{r_1})} + \|\compw\|^2_{\boldL^2_\dens(B_{r_1})})
\end{align}
for some $C_W>0$.
Thus we can combine \eqref{eq:abhierweiter}, \eqref{eq:InequProofG1} and \eqref{eq:InequProofG2} to obtain
\begin{align}\label{eq:EstInt}
\Re\Big( &e^{-i(\theta+\tau)\sign\omega} \spl \op_1^\mathrm{int} \dislag,\dislag \spr_\hs \Big)
\geq C_\mathrm{int} (\|\compv\|_\hs^2
+\|\partial_\bflow\compw\|^2_{\boldL^2_\dens(B_{r_1})} + \|\compw\|^2_{\boldL^2_\dens(B_{r_1})})
\end{align}
with $C_\mathrm{int}:=\cos(\theta+\tau) \min\{C_V-\delta,\,C_W\}>0$.

\emph{coercivity of $A_1$ (2nd part):}
Now we estimate the remaining parts of $A_1$, which are new compared to \cite{HallaHohage:21}.
First we consider
\begin{align*}
\Re \big( e^{-i(\theta+\tau)\sign \omega}
\spl \ol{\sigma}(i\omega\damp+\matii)\compw,\compw'\spr_{A_{r_1,r_2}} \big).
\end{align*}
We recall $\ol{\sigma}=e^{-i\mu}$ and $\matii=\mati+(\omega+i\angvel\times)^*(\omega+i\angvel\times)$.
Since $\mu\geq0$ and $(\omega+i\angvel\times)^*(\omega+i\angvel\times)$ is positive semi-definite and it follows
\begin{align*}
&\sign\omega \arg \big( e^{-i(\theta+\tau)\sign \omega}
\spl \ol{\sigma}(i\omega\damp+\matii)\compw,\compw'\spr_{A_{r_1,r_2}} \big)\\
\leq &
\sign\omega \arg \big( e^{-i(\theta+\tau)\sign \omega}
\spl (i\omega\damp+\mati)\compw,\compw'\spr_{A_{r_1,r_2}} \big)
\leq \sign\omega(\pi/2-\tau)
\end{align*}
from to the definition of $\theta$, see \eqref{eq:theta}.
On the other hand, we can apply $\mu\leq\mu_*$ and use the definition of $\beta$ and \eqref{eq:EstBeta} to estimate
\begin{align}\label{eq:EstArg1}
\begin{split}
&\sign\omega \arg \big( e^{-i(\theta+\tau)\sign \omega}
\spl \ol{\sigma}(i\omega\damp+\matii)\compw,\compw\spr_{B_{r_1}^c} \big)\\
= &
\sign\omega \arg \big( e^{-i(\theta+\tau+\mu)\sign \omega}
\spl (i\omega\damp+\matii)\compw,\compw\spr_{B_{r_1}^c} \big)\\
\geq &
\sign\omega \arg \big( e^{-i(\theta+\tau+\mu_*)\sign \omega}
\spl (i\omega\damp+\matii)\compw,\compw\spr_{B_{r_1}^c} \big)\\
= &
\sign\omega \arg \big( -ie^{-i\beta\sign \omega}
\spl (i|\omega|\damp+\matii)\compw,\compw\spr_{B_{r_1}^c} \big)\\
\geq &
\sign\omega \arg \big( e^{-i\beta\sign\omega} (|\omega|\dampl-i\sign\omega\|M_{\matii}\|_{L(\boldL^2_\dens)} \big)
>-\sign\omega \pi/2.
\end{split}
\end{align}
Hence, there exists a constant $C_1>0$ such that
\begin{align}\label{eq:EstExt1}
\begin{split}
\Re \big( e^{-i(\theta+\tau)\sign \omega}
\spl \ol{\sigma}(i\omega\damp+\matii)\compw,\compw\spr_{A_{r_1,r_2}} \big)
&\geq C_1 |\spl \ol{\sigma}(i\omega\damp+\matii)\compw,\compw\spr_{A_{r_1,r_2}}|\\
&\geq C_1 |\omega|\dampl \spl \compw,\compw\spr_{A_{r_1,r_2}}.
\end{split}
\end{align}
Since also
\begin{align}\label{eq:EstArg2}
\begin{split}
&\sign\omega \arg \big( e^{-i(\theta+\tau)\sign \omega}
\spl \ol{\sigma_*}(i\omega\damp+\matii)\compw,\compw\spr_{B_{r_2}^c} \big)\\
=&\sign\omega \arg \big( e^{-i(\theta+\tau+\mu_*)\sign \omega}
\spl (i\omega\damp+\matii)\compw,\compw\spr_{B_{r_2}^c} \big)\\
=&\sign\omega \arg \big( e^{-i\beta\sign \omega}
\spl (|\omega|\damp-i\matii)\compw,\compw\spr_{B_{r_2}^c} \big)\\
\leq&\sign\omega \arg \big(
|\omega|\dampl+i\sign\omega\|M_{\matii}\|_{L(\boldL^2_\dens)} \big)
<\sign\omega \pi/2,
\end{split}
\end{align}
there exists a constant $C_2>0$ such that
\begin{align}\label{eq:EstExt2}
\begin{split}
\Re \big( e^{-i(\theta+\tau)\sign \omega}
\ol{\sigma_*} \spl (i\omega\damp+\matii)\compw,\compw\spr_{B_{r_2}^c} \big)
&\geq C_2 |\spl (i\omega\damp+\matii)\compw,\compw\spr_{B_{r_2}^c}|\\
&\geq C_2 |\omega|\dampl \spl \compw,\compw\spr_{B_{r_2}^c}.
\end{split}
\end{align}
Further, we compute
\begin{align*}
&\Re \big( -e^{-i(\theta+\tau)\sign \omega}\ol{\sigma_*} \big)
= \Re \big( i\sign\omega e^{-i\sign\omega \beta}\big)
= \sin\beta
\end{align*}
and hence
\begin{align}\label{eq:EstExt3}
\begin{split}
\Re \big( -e^{-i(\theta+\tau)\sign \omega}
\ol{\sigma_*} \spl \sounds^2(\div+\qcoeff\cdot)\compw,(\div+\qcoeff\cdot)\compw\spr_{B_{r_2}^c}
\big) &=\sin\beta
\| \sounds(\div+\qcoeff\cdot)\compw\|_{L^2_\dens({B_{r_2}^c})}^2,\\
\Re \big( -e^{-i(\theta+\tau)\sign \omega}
\ol{\sigma_*}
\spl \compz,\compz \spr
\big) &=\sin\beta
\|\compz\|^2_\hs.
\end{split}
\end{align}
We combine \eqref{eq:EstInt}, \eqref{eq:EstExt1}, \eqref{eq:EstExt2} and \eqref{eq:EstExt3} to obtain
\begin{align}\label{eq:Total1}
\begin{split}
\Re\Big( &e^{-i(\theta+\tau)\sign\omega} \spl \op_1 \dislag,\dislag \spr_\hs \Big)
\geq\\
& \hspace{4mm} C_3 (\|\compv\|_\hs^2
+\|\partial_\bflow\compw\|^2_{\boldL^2_\dens}
+\|\compw\|^2_{\boldL^2_\dens}
+\|(\div+\qcoeff\cdot)\compw\|_{L^2_\dens({B_{r_2}^c})}^2
+\|\compz\|^2_\hs)
\end{split}
\end{align}
with $C_3:=\min\{C_\mathrm{int},C_1|\omega|\dampl,C_2|\omega|\dampl,\sin\beta,\soundsl^2\sin\beta\}>0$.
In $B_{r_2}$ it holds for $\compw\in\hsW$ that $\div\compw=-\qcoeff\cdot\compw$ and hence
\begin{align}\label{eq:EstDiv1}
\|\div\compw\|_{L^2_\dens(B_{r_2})}
=\|\qcoeff\cdot\compw\|_{L^2_\dens(B_{r_2})}
\leq \max\{1,\|\qcoeff\|_{\boldL^\infty}\} \|\compw\|_{L^2_\dens(B_{r_2})}.
\end{align}
On the other hand a weighted Cauchy-Schwarz inequality yields the existence of a constant $C_4>0$ such that
\begin{align}\label{eq:EstDiv2}
\|(\div+\qcoeff\cdot)\compw\|_{L^2_\dens(B_{r_2}^c)}^2 +\|\compw\|_{L^2_\dens(B_{r_2}^c)}^2
\geq C_4 (\|\div\compw\|_{L^2_\dens(B_{r_2}^c)}^2 +\|\compw\|_{L^2_\dens(B_{r_2}^c)}^2).
\end{align}
Finally, we combine \eqref{eq:Total1}, \eqref{eq:EstDiv1} and \eqref{eq:EstDiv2} to
\begin{align*}
\Re\Big( &e^{-i(\theta+\tau)\sign\omega} \spl \op_1 \dislag,\dislag \spr_\hs \Big)
\geq\\
& \hspace{4mm} C_5 (\|\compv\|_\hs^2
+\|\partial_\bflow\compw\|^2_{\boldL^2_\dens}
+\|\compw\|^2_{\boldL^2_\dens}
+\|\div\compw\|_{L^2_\dens}^2
+\|\compz\|^2_\hs)\\
&=C_5 (\|\compv\|_\hs^2+\|\compw\|^2_\hs+\|\compz\|^2_\hs)
\end{align*}
with $C_5:=C_3 \min\{C_4, \,1/\max\{1,\|\qcoeff\|_{\boldL^\infty}\}\}/2>0$.
Since the decomposition $\hs=\hsV\oplus^\calT\hsW\oplus^\calT\hsZ$ is a topological one, it follows that $A_1$ is coercive.
\end{proof}

\begin{proposition}\label{prop:bij-cow}
Let the assumptions of Theorem~\ref{thm:wTc-cow} be satisfied.
Then $\opAcow$ is bijective.
\end{proposition}
\begin{proof}
Follows from Proposition~\ref{prop:Fredholm}, Lemma~\ref{lem:inj-cow} and Theorem~\ref{thm:wTc-cow}.
\end{proof}

\begin{remark}
In Theorem~\ref{thm:wTc-cow} the radius $r_2$ can be chosen arbitrarily close to $r_1$.
Further, the balls $B_{r_1}$, $B_{r_2}$ can be replaced by simply connected open Lipschitz domains $D_1, D_2$ whereby $D_1$ contains $\supp\bflow$, $D_2$ contains the closure of $D_1$ and $D_2$ is of class $C^{1,1}$.
Thus in the definition \eqref{eq:theta} of $\theta$ the ball $B_{r_2}$ can be replaced by a smooth, simply connected, arbitrarily close neighborhood of $\supp\bflow$.
\end{remark}

\section{The full equation}\label{sec:full}

Let us now consider the full equation \eqref{eq:PDE}, which includes the Eulerian perturbation $\graveu$ of the gravitational background potential $\grav$.
First we note the injectivity of $\op$.
\begin{lemma}\label{lem:inj-full}
Let the assumptions of Section~\ref{subsec:assumptions} be satisfied and $\omega\neq0$.
Then $\op$ is injective.
\end{lemma}
\begin{proof}
Let $(\dislag,\graveu)\in\ker A$. Then 
\begin{align*}
0=\big|\Im \big(\ses\big((\dislag,\graveu),(\dislag,\graveu)\big)\big)\big|
= |\omega| \spl \damp \dislag,\dislag \spr
\geq |\omega|\dampl \|\dislag\|^2_{\boldL^2_\dens}\end{align*}
and hence $\dislag=0$. We further compute 
\begin{align*}
0=\ses\big((\dislag,\graveu),(\dislag,\graveu)\big)
= \ses\big((0,\graveu),(0,\graveu)\big)
= \frac{1}{4\pi G}\|\graveu\|^2_{\hsgrav}
\end{align*}
and conclude that $\graveu=0$. Hence, $(\dislag,\graveu)=(0,0)$.
\end{proof}
Since $a\big((0,\graveu),(0,\graveu')\big)=\frac{1}{4\pi G}\spl \graveu,\graveu' \spr_{\hsgrav}$, we can build the Schur complement $A_\mathrm{Schur}$ of $A$ with respect to $\graveu$.
Thus $A$ is Fredholm if and only if $A_\mathrm{Schur}$ is so, and hence it suffices to analyze $A_\mathrm{Schur}$.
The difference $\opAcow-A_\mathrm{Schur}$ is of the form $E^*RE$ with the embedding $E\colon\hs\rightarrow\boldL^2_\dens$ and a positive semi-definite operator $R\in L(\boldL^2_\dens)$.
Thus, to reuse our analysis for $\opAcow$ we can put the term $\spl RE\boldw,E\boldw'\spr_\hs$
together with $\spl(i\damp+\mati)\boldw,\boldw'\spr$.
In particular, let $M_{i\damp+\mati}\in L(\boldL^2_\dens)$ be the multiplication operator with symbol $i\damp+\mati$.
The we can get control over the numerical range of $E^*(M_{i\damp+\mati}+R)E$.
If we continue this approach we need to consider the form of our test function $\dislag'=T\dislag$, and consequently we have to treat the numerical range of the term $\spl (M_{i\damp+\mati}+R)\boldw,\sigma\boldw\spr$.
The difficulty here is that $R$ is a non local operator and hence its interaction with the multiplication by $\sigma$ is unclear.
To bypass this obstacle we split of some compact operators such that the embedding $E\colon\hs\rightarrow\boldL^2_\dens$ can be replaced by the embedding $E_2\colon\hs\rightarrow\boldL^2_\dens(B_{r_2}^c)$, and $R\in L(\boldL^2_\dens)$ is replaced by some other positive semi-definite operator $\tilde R\in L(\boldL^2_\dens(B_{r_2}^c))$.
Since $\sigma$ is constant in $B_{r_2}^c$, it is then easy to control the numerical range of
$\spl (M_{i\damp+\mati}+\tilde R)\boldw,\sigma\boldw\spr_{B_{r_2}^c}
=\ol{\sigma_*} \spl (M_{i\damp+\mati}+\tilde R)\boldw,\boldw\spr_{B_{r_2}^c}$.

In the following we carry out this approach in detail.
To this end we assume that there exists $r_3>r_2>r_1$ such that $\dens\in W^{1,\infty}(B_{r_3})$.
Then let
\begin{align*}
H_1 &:= \{\graveu\in \hsgrav\colon \graveu=0 \text{ in } B_{r_2}\},\qquad
H_2 := H_1^{\bot_{\hsgrav}}\subset\{\graveu\in \hsgrav\colon \Delta\graveu=0 \text{ in } B_{r_2}^c\},
\end{align*}
with associated orthogonal projections $\PHi$ and $\PHii$.
Let $Q\in L(\hsgrav,\boldL^2_\dens(B_{r_2}^c))$ be defined by
\begin{align*}
\spl \boldxi,Q\graveu' \spr_{B_{r_2}^c}&:=-\spl \dens\boldxi,\nabla \PHi\graveu' \spr_{\boldL^2(B_{r_2}^c)}
-\spl \dens\boldxi,\nabla \PHii\graveu' \spr_{\boldL^2(B_{r_3}^c)}
\end{align*}
for all $\boldxi\in \boldL^2_\dens(B_{r_2}^c)$, $\graveu'\in\hsgrav$,
and $\Kf\in L(\hsgrav,\hs)$ be defined by
\begin{align*}
\spl \dislag,\Kf\graveu' \spr_\hs&:=\spl (\dens\div+\nabla\dens\,\cdot)\dislag, \PHii\graveu' \spr_{\boldL^2(B_{r_3})}
-\spl \dens\,\nv\cdot\dislag, \PHii\graveu' \spr_{H^{-1/2}(\partial B_{r_3})\times H^{1/2}(\partial B_{r_3})},
\end{align*}
for all $\dislag\in\hs$, $\graveu'\in\hsgrav$.
Further, let $E_2\in L(\hs, \boldL^2_\dens(B_{r_2}^c))$ be the respective embedding operator.
Then by means of $Q, K$ and $E_2$ we can express
\begin{align*}
a\big((\dislag,0),(0,\graveu')\big)
&=-\spl \dislag,\nabla\graveu' \spr
=-\spl \dens\dislag,\nabla\graveu' \spr_{\boldL^2}\\
&=-\spl \dens\dislag,\nabla \PHi\graveu' \spr_{\boldL^2}
-\spl \dens\dislag,\nabla \PHii\graveu' \spr_{\boldL^2}\\
&=-\spl \dens\dislag,\nabla \PHi\graveu' \spr_{\boldL^2(B_{r_2}^c)}
-\spl \dens\dislag,\nabla \PHii\graveu' \spr_{\boldL^2(B_{r_3}^c)}
-\spl \dens\dislag,\nabla \PHii\graveu' \spr_{\boldL^2(B_{r_3})}\\
&=-\spl \dens\dislag,\nabla \PHi\graveu' \spr_{\boldL^2(B_{r_2}^c)}
-\spl \dens\dislag,\nabla \PHii\graveu' \spr_{\boldL^2(B_{r_3}^c)}\\
&+\spl (\dens\div+\nabla\dens\,\cdot)\dislag, \PHii\graveu' \spr_{\boldL^2(B_{r_3})}
-\spl \dens\,\nv\cdot\dislag, \PHii\graveu' \spr_{H^{-1/2}(\partial B_{r_3})\times H^{1/2}(\partial B_{r_3})}\\
&=\spl E_2\dislag,Q\graveu' \spr_{B_{r_2}^c}
+\spl \dislag,\Kf\graveu' \spr_\hs\\
&=\spl \dislag,E_2^*Q\graveu' \spr_\hs
+\spl \dislag,\Kf\graveu' \spr_\hs
\end{align*}
and like-wise $a\big((0,\graveu),(\dislag',0)\big)=\spl E_2^*Q\graveu,\dislag' \spr_\hs
+\spl \Kf\graveu,\dislag' \spr_\hs$.
Due to the compact Sobolev embedding $H^1(B_{r_3})\hookrightarrow L^2(B_{r_3})$ the operator related to the first term $\spl (\dens\div+\nabla\dens\,\cdot)\dislag, \PHii\graveu' \spr_{\boldL^2(B_{r_3})}$ in the definition of $\Kf$ is compact.
Since $\PHii\graveu'$ solves $\Delta \PHii\graveu'=0$ in $B_{r_2}^c$ and $r_3>r_2$, it follows with standard regularity theory (e.g.\ \cite[Theorems 2.3.3.2, 2.4.2.7]{Grisvard:85} applied to $\chi\graveu$ with a suitable cut-off function $\chi$) that the embedding $H_2\to H^2(A_{(r_2+r_3)/2,r_3})$ is bounded.
Hence, the operator related to the second term $-\spl \dens\,\nv\cdot\dislag, \PHii\graveu' \spr_{H^{-1/2}(\partial B_{r_3})\times H^{1/2}(\partial B_{r_3})}$ in the definition of $\Kf$ is compact too.
Together, it follows that $\Kf$ is compact.
Since the Fredholmness and the index of an operator are invariant under compact perturbations, it suffices to analyze $a\big((\dislag,\graveu),(\dislag',\graveu')\big)-\spl \dislag,\Kf\graveu' \spr_\hs-\spl \Kf\graveu,\dislag' \spr_\hs$ instead of $a\big((\dislag,\graveu),(\dislag',\graveu')\big)$.
As $a\big((0,\graveu),(0,\graveu')\big)=\frac{1}{4\pi G}\spl \graveu,\graveu' \spr_{\hsgrav}$ is unchanged under the previous perturbation, we can again build the Schur complement of $a\big((\dislag,\graveu),(\dislag',\graveu')\big)-\spl \dislag,\Kf\graveu' \spr_\hs-\spl \Kf\graveu,\dislag' \spr_\hs$ with respect to $\graveu$ and obtain
\begin{align*}
A_\mathrm{Schur}=\opAcow-4\pi G \, E_2^*QQ^*E_2.
\end{align*}

\begin{lemma}\label{lem:wTc-schur}
Let the assumptions of Section~\ref{subsec:assumptions} be satisfied.
Let $r_3>r_2>r_1$ be such that $\sounds, \dens \in W^{1,\infty}(B_{r_3})$ and $\omega\neq0$.
If
\begin{align*}
\|\sounds^{-1}\bflow\|_{\boldL^\infty}^2<\frac{1}{1+\tan^2\theta},
\end{align*}
then there exists $\mu$ such that $\opAcow-4\pi G \, E_2^*QQ^*E_2$ is weakly $T$-coercive.
\end{lemma}
\begin{proof}
To prove that $\opAcow-4\pi G \, E_2^*QQ^*E_2$ is weakly $T$-coercive we can apply the same technique as in the proof of Theorem~\ref{thm:wTc-cow} and only require to make some minor adaptations.
We need to split $T^*(\opAcow-4\pi G \, E_2^*QQ^*E_2)=\tilde A_1+\tilde A_2$ into a coercive operator $\tilde A_1$ and a compact operator $\tilde A_2$.
Let $T^*\opAcow=A_1+A_2$ be the decomposition as in the proof of Theorem~\ref{thm:wTc-cow}.
Then let
\begin{align*}
\tilde A_1&:=A_1+4\pi G \ol{\sigma_*} \PW^*E_2^*QQ^*E_2\PW,\\
\tilde A_2&:=A_2-4\pi G \big(
\PZ^*E_2^*QQ^*E_2\PW
-\ol{\sigma_*} \PW^*E_2^*QQ^*E_2\PZ
+\PZ^*E_2^*QQ^*E_2\PZ
\big).
\end{align*}
Note that since $\boldv=0$ in $B_{r_2}^c$ there arise in $\tilde A_1$ and $\tilde A_2$ compared to $A_1$ and $A_2$ no additional terms involving $\boldv$. 
Further, $\tilde A_2$ is indeed compact, because of the compactness of $A_2$ and $\PZ$.
The additional operator $4\pi G \ol{\sigma_*} \PW^*E_2^*QQ^*E_2\PW$ in $\tilde A_1$ can be estimated together with $\ol{\sigma_*}\PW^*E_2^*M_{i\damp+\matii}E_2\PW$.
To this end $\beta$ needs to be chosen such that
\begin{align}\label{eq:EstBetaFull}
\begin{split}
0&<
\Re \big( -ie^{-i\sign\omega\beta} (i|\omega|\dampl+\sign\omega\|M_{\matii}+QQ^*\|_{L(\boldL^2_\dens(B_{r_2}^c))} \big)\\
&=\Re \big( e^{-i\beta} (|\omega|\dampl-i\|M_{\matii}+QQ^*\|_{L(\boldL^2_\dens(B_{r_2}^c))} \big).
\end{split}
\end{align}
is satisfied instead of \eqref{eq:EstBeta}.
Then in addition to \eqref{eq:EstArg1} and \eqref{eq:EstArg2} we apply \eqref{eq:EstBetaFull} to estimate
\begin{align*}
&\sign\omega \arg \big( e^{-i(\theta+\tau)\sign \omega}
\ol{\sigma_*}
\spl (i\omega M_\damp+M_{\matii}+QQ^*)\compw,\compw\spr_{B_{r_2}^c} \big)\\
\geq &
\sign\omega \arg \big( e^{-i\beta\sign\omega} (|\omega|\dampl-i\sign\omega\|M_{\matii}+QQ^*\|_{L(\boldL^2_\dens(B_{r_2}^c))} \big)
>-\sign\omega \pi/2.
\end{align*}
and
\begin{align*}
&\sign\omega \arg \big( e^{-i(\theta+\tau)\sign \omega}
\ol{\sigma_*}
\spl (i\omega M_\damp+M_{\matii}+QQ^*)\compw,\compw\spr_{B_{r_2}^c} \big)\\
\leq&\sign\omega \arg \big(
|\omega|\dampl+i\sign\omega\|M_{\matii}+QQ^*\|_{L(\boldL^2_\dens(B_{r_2}^c))} \big)
<\sign\omega \pi/2.
\end{align*}
The remainder of the proof is identical to the proof of Theorem~\ref{thm:wTc-cow}.
\end{proof}

We can now prove the first main result of this article:

\begin{theorem}\label{thm:bij-full}
Let the parameters $\omega, \angvel, \sounds, \bflow, \dens, \pres, \grav$ and $\damp$ satisfy the assumptions of Section~\ref{subsec:assumptions}, in particular $\supp\bflow\subset B_{r_1}$.
Let $r_3>r_2>r_1$ be such that $\sounds, \dens \in W^{1,\infty}(B_{r_3})$, $\omega\neq0$ and $\theta$ be as defined in \eqref{eq:theta}.
Let the Hilbert spaces $\hs$ and $\hsgrav$ be as defined in \eqref{eq:hsX} and \eqref{eq:hsgrav} respectively, let the sesquilinearform $\ses(\cdot,\cdot)$ be as defined in \eqref{eq:sesq-a} and $\op\in\blo(\hs\times\hsgrav)$ be the associated operator as defined in \eqref{eq:defRieszRep}.
If
\begin{align*}
\|\sounds^{-1}\bflow\|_{\boldL^\infty}^2<\frac{1}{1+\tan^2\theta},
\end{align*}
then $\op$ is bijective, and hence for each $\boldf\in\boldL^2_\dens$ Equation \eqref{eq:PDE} admits a unique solution $(\dislag,\graveu)\in \hs\times\hsgrav$ which depends continuously on $\boldf$.
\end{theorem}
\begin{proof}
Recall that $A$ is injective due to Lemma~\ref{lem:inj-full} and hence it suffices to show that $A$ is Fredholm with index zero.
Since the Fredholmness and the index of an operator are unchanged by compact perturbations and because $K_3$ is compact, $A$ is Fredholm with index zero if and only if $A-\bpm 0&\Kf\\ \Kf^*&0\epm$ is so.
The second diagonal component of $A-\bpm 0&\Kf\\ \Kf^*&0\epm$ equals $\frac{1}{4\pi G} I_{\hsgrav}$ and hence we build the Schur complement of $A-\bpm 0&\Kf\\ \Kf^*&0\epm$ with respect to $\graveu$.
Thus $A$ is Fredholm with index zero, if and only if the Schur complement $\opAcow-4\pi G \, E_2^*QQ^*E_2$ is so.
However, $\opAcow-4\pi G \, E_2^*QQ^*E_2$ is weakly $T$-coercive due to Lemma~\ref{lem:wTc-schur}.
Thus the claim is proven.
\end{proof}

\section{A scalar equation in the atmosphere}\label{sec:cp}

Consider spherical variables $r=|\boldx|$ and $\hat\boldx=|\boldx|^{-1}\boldx$.
Assume that $\supp\boldf\subset B_{r_2}$ and that in $B_{r_2}^c$ the parameters $\sounds$, $\dens$ and $\pres$ depend only on $r$ and let us use the sloppy notation $\sounds(r)=\sounds(\boldx)$, etc..
For the following discussion we consider only the exterior domain $B_{r_2}^c$.
Under the previous assumptions it holds
\begin{align*}
\qcoeff(\boldx)=\frac{\partial_r \pres(r)}{\sounds^2(r)\dens(r)} \hat\boldx.
\end{align*}
Let
\begin{align*}
\eta(r):=\int_{r_2}^r \frac{\partial_r \pres(r')}{\sounds^2(r')\dens(r')} \dd r'.
\end{align*}
Since $\qcoeff\in\boldL^\infty$ it follows $\frac{\partial_r \pres}{\sounds^2\dens}\in L^\infty(r_2,\infty)$ and hence $\eta$ is well defined and $\eta\in C([r_2,\infty))$.
Then it follows $\qcoeff=\nabla\eta$ and we can obtain the representation
\begin{align*}
(\div+\qcoeff\cdot)\boldu=e^{-\eta}\div(e^\eta \boldu).
\end{align*}
Hence in $B_{r_2}^c$ the equation for the Cowling approximation reads
\begin{align}\label{eq:EqinAtmo}
-e^\eta \nabla (\sounds^2\dens e^{-2\eta} \div (e^\eta \boldu))
-\dens(\matii+i\omega\damp) \boldu =0.
\end{align}
Thus $\boldu$ satisfies
\begin{align*}
\curl(e^{-\eta}\dens(\matii+i\omega\damp) \boldu)=0.
\end{align*}
In the following we assume $\omega\neq0$.
Then for each $\boldx\in B_{r_2}^c$ the matrix $\matii(\boldx)+i\omega\damp(\boldx) I_{3\times3}$ is coercive and hence invertible.
Thus it follows that
\begin{align}\label{eq:UeqGrad}
\boldu=e^{\eta}\dens^{-1} (\matii+i\omega\damp)^{-1} \nabla v
\end{align}
with a scalar function $v$.
We plug \eqref{eq:UeqGrad} into \eqref{eq:EqinAtmo} and obtain
\begin{align*}
-e^\eta \nabla \Big( \sounds^2\dens e^{-2\eta} \div (e^{2\eta}\dens^{-1} (\matii+i\omega\damp)^{-1} \nabla v) +v\Big)=0
\end{align*}
and hence $v$ satisfies
\begin{align*}
\sounds^2\dens e^{-2\eta} \div (e^{2\eta}\dens^{-1} (\matii+i\omega\damp)^{-1} \nabla v) +v=\mathrm{const}.
\end{align*}
If $\mathrm{const}\neq0$ the function $v\equiv\mathrm{const}$ solves the former equation.
However, for $v\equiv\mathrm{const}$ it holds $\boldu=e^{\eta}\dens^{-1} (\matii+i\omega\damp)^{-1} \nabla v=0$ and thus it is sufficient to consider only the homogeneous equation.
A reformulation yields
\begin{align*}
-\div \left(\frac{e^{2\eta}}{\dens} (\matii+i\omega\damp)^{-1} \nabla v\right)
-\frac{e^{2\eta}}{\sounds^2\dens} v=0.
\end{align*}
So we couple the equation for $\boldu$ in $B_{r_2}$ with the equation for $v$ in $B_{r_2}^c$ and end up with the system
\begin{subequations}
\label{eq:PDE-cp}
\begin{align}
\begin{split}
\label{eq:PDEint}
-\nabla \big(\sounds^2 \dens\div \dislag
+\nabla \pres\cdot\dislag\big)
+\nabla \pres \div \dislag
+\hess(\pres)\dislag
& \phantom{= \dens\boldf \quad\text{in }\setR^3,}\\
-\dens\hess(\grav)\dislag
-\dens\big(\omega +i\partial_{\bflow}+i\angvel\times \big)^2 \dislag
-i\omega\damp\dens\dislag
&= \dens\boldf \quad\text{in } B_{r_2},
\end{split}\\
\label{eq:PDEext}
-\div \left(\frac{e^{2\eta}}{\dens} (\matii+i\omega\damp)^{-1} \nabla v\right)
-\frac{e^{2\eta}}{\sounds^2\dens} v&=\phantom{\dens}0 \quad\text{in } B_{r_2}^c,\\
\label{eq:IntF1}
\nv\cdot e^\eta\dens^{-1} (\matii+i\omega\damp)^{-1} \nabla v&=\nv\cdot\boldu \quad\text{on } \partial B_{r_2},\\
\label{eq:IntF2}
\div\left(e^\eta\dens^{-1} (\matii+i\omega\damp)^{-1} \nabla v\right)&=\div\boldu \quad\text{on } \partial B_{r_2}.
\end{align}
\end{subequations}
We introduce the sesquilinearforms
\begin{align*}
\begin{split}
\sescow^\mathrm{int}(\dislag,\dislag')
&:=\spl \sounds^2 (\div \dislag + \qcoeff\cdot\dislag), \div \dislag'  + \qcoeff\cdot\dislag'\spr_{B_{r_2}}\\
&-\spl (\omega+i\partial_\bflow+i\angvel\times) \dislag, (\omega+i\partial_\bflow+i\angvel\times) \dislag' \spr_{B_{r_2}}\\
&-\spl \mati \dislag,\dislag' \spr_{B_{r_2}}
-i\omega \spl \damp \dislag, \dislag' \spr_{B_{r_2}},
\end{split}
\end{align*}
and
\begin{align}\label{eq:acp}
\begin{split}
\sescp\big((\dislag,v),(\dislag',v')\big)&:=
\sescow^\mathrm{int}(\dislag,\dislag')\\
&+\spl \nv\cdot\boldu,v'\spr_{H^{-1/2}(\partial B_{r_2})\times H^{1/2}(\partial B_{r_2})}\\
&+\spl v,\nv\cdot\boldu'\spr_{H^{1/2}(\partial B_{r_2})\times H^{-1/2}(\partial B_{r_2})}\\
&+\Big\spl \frac{e^{2\eta}}{\dens} (\matii+i\omega\damp)^{-1} \nabla v,\nabla v'\Big\spr_{\boldL^2(B_{r_2}^c)}
-\Big\spl \frac{e^{2\eta}}{\sounds^2\dens} v,v' \Big\spr_{\boldL^2(B_{r_2}^c)}.
\end{split}
\end{align}
Further, we introduce the Hilbert spaces
\begin{subequations}\label{eq:DefXintVext}
\begin{align}
\Xint&:=\hs(B_{r_2}),\\
\begin{split}
\Vext&:=\{v \in L_{e^{2\eta}/\dens}^2(B_{r_2}^c)\colon \quad \nabla v \in \boldL_{e^{2\eta}/\dens}^2(B_{r_2}^c)\},\\
\spl v,v' \spr_{\Vext} &:= \spl v, v' \spr_{L_{e^{2\eta}/\dens}^2(B_{r_2}^c)}
+\spl \nabla v, \nabla v' \spr_{\boldL_{e^{2\eta}/\dens}^2(B_{r_2}^c)}.
\end{split}
\end{align}
\end{subequations}
Then the weak formulation of \eqref{eq:PDE-cp} is to find $(\boldu,v)\in \Xint\times\Vext$ such that
\begin{align*}
\sescp\big((\dislag,v),(\dislag',v')\big)=\spl \boldf,\boldu' \spr_{B_{r_2}}
\end{align*}
for all $(\boldu',v')\in \Xint\times\Vext$.
In particular, if we test the left-hand-side of \eqref{eq:PDEext} with $v'\in\Vext$ and integrate by parts we obtain the last line of the right-hand-side of \eqref{eq:acp} plus the boundary term
\begin{align*}
-\Big\spl \nv\cdot \frac{e^{2\eta}}{\dens} (\matii+i\omega\damp)^{-1} \nabla v, v'\Big\spr_{H^{-1/2}(\partial B_{r_2}^c) \times H^{1/2}(\partial B_{r_2}^c)}.
\end{align*}
Note that for the outward unit normal vectors $\nv$ of $\partial B_{r_2}$ and $\nv'$ of $\partial B_{r_2}^c$ it holds $\nv'=-\nv$.
Hence with \eqref{eq:IntF1} and $\eta(r_2)=0$ the boundary term equals the second line in the right-hand-side of \eqref{eq:acp}.
If we test the left-hand-side of \eqref{eq:PDEint} with $\dislag'\in\Xint$ and integrate by parts we obtain the first line of the right-hand-side of \eqref{eq:acp} plus the boundary term
\begin{align*}
- \Big\spl \sounds^2\dens (\div+\qcoeff\cdot)\dislag, \nv\cdot\dislag'\Big\spr_{H^{1/2}(\partial B_{r_2}^c) \times H^{-1/2}(\partial B_{r_2}^c)}.
\end{align*}
Since on $\partial B_{r_2}$ the vectors $\qcoeff$ and $\nv$ are parallel we obtain on $\partial B_{r_2}$ with \eqref{eq:IntF1} and \eqref{eq:IntF2} that
\begin{align*}
\dens\sounds^2(\div+\qcoeff\cdot)\boldu
&=\dens\sounds^2(\div+\qcoeff\cdot) e^{\eta}\dens^{-1} (\matii+i\omega\damp)^{-1} \nabla v\\
&=\dens\sounds^2e^{-\eta}\div(e^{2\eta}\dens^{-1} (\matii+i\omega\damp)^{-1} \nabla v)
=-e^\eta v
=-v,
\end{align*}
whereat we used \eqref{eq:PDEext} for the third equality and $\eta(r_2)=0$ for the last equality.
Hence the boundary term equals the third line of the right-hand-side of \eqref{eq:acp}.
The injectivity of $\Acp$ can be seen similarly as for $A$ and $\opAcow$.
\begin{lemma}\label{lem:inj-cp}
Let the assumptions of Section~\ref{subsec:assumptions} be satisfied and $\omega\neq0$.
Then $\Acp$ is injective.
\end{lemma}
\begin{proof}
Let $(\dislag,v)\in\ker\Acp$. Then 
\begin{align*}
\Im \big(\sescp((\dislag,v),(\dislag,v))\big)
=\Im\Big(\Big\spl \frac{e^{2\eta}}{\dens} (\matii+i\omega\damp)^{-1} \nabla v,\nabla v\Big\spr_{\boldL^2(B_{r_2}^c)}\Big)
-\omega \spl \damp \dislag,\dislag \spr_{B_{r_1}}.
\end{align*}
Since $\matii(\boldx)$ is a selfadjoint matrix for each $\boldx\in B_{r_2}$ and $\damp\geq\dampl>0$, it follows that
\begin{align*}
0=|\Im \big(\sescp((\dislag,v),(\dislag,v))\big)|
\geq |\omega|\dampl \spl \dislag,\dislag \spr_{B_{r_1}}
+ C \Big\spl \frac{e^{2\eta}}{\dens} \nabla v,\nabla v\Big\spr_{\boldL^2(B_{r_2}^c)}
\end{align*}
with a constant $C>0$.
Thus $\dislag=0$ and $v\equiv v_c\in\setC$ is constant.
Hence
\begin{align*}
0=\sescp((\dislag,v),(\dislag,v))
=\sescp((0,v_c),(0,v_c))
=-\Big\spl \frac{e^{2\eta}}{\sounds^2\dens} v_c,v_c \Big\spr_{\boldL^2(B_{r_2}^c)}
\end{align*}
and thus $v_c=0$.
Altogether it follows $(\dislag,v)=0$ and the claim is proven.
\end{proof}
Next we introduce a topological decomposition of $\Xint$ in analogy to Theorem~\ref{thm:decomposition}.
\begin{theorem}\label{thm:decomposition-Xint}
Let $\dens$, $\bflow$, $\qcoeff$ and $r_1$ be as in Section~\ref{subsec:assumptions} and $r_2>r_1$.
Then $\Xint=\hs(B_{r_2})$ admits a topological decomposition  
\begin{align*}
\Xint = \hsVint \oplus^\calT \hsWint \oplus^\calT \hsZint
\end{align*}
with the following properties:
\begin{enumerate}
\item\label{it:VcompactembeddingInt}
$\hsVint\subset\{\nabla\compvp\colon
\compvp\in H^2(B_{r_2})\mbox{ with }\frac{\partial \compvp}{\partial \nv}=0 \mbox{ on }\partial B_{r_2}\}$ is compactly embedded in $\boldL^2(B_{r_2})$.
\item\label{it:WkernelInt}
$\hsWint = \{\dislag \in \Xint \colon \div \dislag + \qcoeff\cdot \dislag = 0 \text{ in } B_{r_2} \} $.
\item\label{it:ZfiniteInt}
$\hsZint$ is finite-dimensional. 
\end{enumerate}
Moreover, for each $\zeta\in W^{1,\infty}(B_{r_2})$ there exists a compact operator $\tilde K_\zeta\in\blo(\Xint)$ such that
\begin{align*}
\|\zeta\div\compv\|_{L^2(B_{r_2})}^2=
\|\zeta\nabla\compv\|_{(L^2(B_{r_2}))^{3x3}}^2
+\spl \tilde K_\zeta\compv,\compv \spr_{\hs}
\end{align*}
for all $\compv\in\hsV$.
\end{theorem}
\begin{proof}
Follows as for Theorem~\ref{thm:decomposition} and \cite[Theorem~3.5]{HallaHohage:21}.
\end{proof}
To construct a suitable T-operator let $r_3>r_2$ and $\mu\in C^\infty(\setR)$ with
\begin{enumerate}
 \item $\mu(r)=0$ for $r\leq r_1$,
 \item $\mu$ is non-decreasing in $A_{r_1,r_2}$,
 \item $\mu(r_2)\in(0,\pi)$,
 \item $\mu$ is non-increasing in $A_{r_2,r_3}$,
 \item $\mu(r)=\mu_*$ for $r\geq r_3$ with constant $\mu_*\in(0,\pi)$.
\end{enumerate}
For such $\mu$ let $\sigma(\boldx):=e^{i\mu(|\boldx|)\sign\omega}$ and $\sigma_*:=e^{i\mu_*\sign\omega}$.
Then we construct $\Tint\in L(\Xint)$ in analogy to $T$ in Section~\ref{sec:cowling}.
In particular, for given $\compw\in\hsWint$ we consider the problem to find $\check v_0\in H^1(A_{r_1,r_2})$ such that
\begin{align*}
(\div+\qcoeff\cdot)(\nabla \check v_0-\sigma\compw)&=0 \quad\text{in } A_{r_1,r_2},\\
\nv\cdot\nabla \check v_0&=0 \quad\text{on }\partial A_{r_1,r_2}.
\end{align*}
Since this equation is weakly coercive, we can find a finite dimensional subspace $\hsWint^0\subset\hsWint$ and a projection $\tilde K_1\in L(\hsWint,\hsWint^0)$ onto $\hsWint^0$ such that the problem to find $\check v_0\in H^1(A_{r_1,r_2})$ with
\begin{align*}
(\div+\qcoeff\cdot)(\nabla \check v_0-\sigma\compw+\sigma \tilde K_1\compw)&=0 \quad\text{in } A_{r_1,r_2},\\
\nv\cdot\nabla \check v_0&=0 \quad\text{on }\partial A_{r_1,r_2},
\end{align*}
admits a unique solution $\check v_0$.
With convenient regularity theory (e.g.\ \cite[Theorems 2.3.3.2, 2.4.2.7]{Grisvard:85}) it follows that $\check v_0 \in H^2(A_{r_1,r_2})$ and $\|\check v_0\|_{H^2(A_{r_1,r_2})} \lesssim \|\boldu\|_{\Xint}$.
Let $\check \boldv$ be the continuation of $\nabla \check v_0$ to $B_{r_1}$ by zero.
It follows that $\check \boldv\in\Xint$.
Hence let
\begin{align*}
T_{\hsWint}'\boldw:=\sigma\boldw-\check\boldv-\sigma \tilde K_1\boldw.
\end{align*}
From its definition it follows that $(\div+\qcoeff\cdot)T_{\hsWint}'\boldw=0$ in $A_{r_1,r_2}$.
Since $\sigma=1$ in $B_{r_1}$, $\tilde K_1$ maps into $\hsWint^0\subset\hsWint$ and $\check\boldv=0$ in $B_{r_1}$ it follows $(\div+\qcoeff\cdot)T_{\hsWint}'\boldw=0$ in $B_{r_1}$ too.
Hence $T_{\hsWint}'\in L(\hsWint)$.
We compute
\begin{align*}
\spl T_{\hsWint}'\boldw,\boldw \spr_\hs&=
\|\boldw\|_{\hs(B_{r_1})}^2
-\spl \tilde K_1\boldw,\boldw \spr_{\hs(B_{r_1})}\\
&+\spl \sigma\boldw,\boldw \spr_{A_{r_1,r_2}}
-\spl \check\boldv,\boldw \spr_{A_{r_1,r_2}}
-\spl \sigma \tilde K_1\boldw,\boldw \spr_{A_{r_1,r_2}}\\
&+\spl \sigma\qcoeff\cdot\boldw,\qcoeff\cdot\boldw \spr_{A_{r_1,r_2}}
-\spl \qcoeff\cdot\check\boldv,\qcoeff\cdot\boldw \spr_{A_{r_1,r_2}}
-\spl \qcoeff\cdot \sigma \tilde K_1\boldw,\qcoeff\cdot\boldw \spr_{A_{r_1,r_2}}
\end{align*}
All terms which involve $\tilde K_1$ are compact, because $\tilde K_1$ is compact.
Also, all terms which involve $\check\boldv$ are compact due to the compact embedding $\boldH^1(A_{r_1,r_2})\hookrightarrow \boldL^2(A_{r_1,r_2})$.
With  $\div\boldw=-\qcoeff\cdot\boldw$ in $B_{r_2}$ for $\boldw\in\hsW$, we obtain that the remainder
\begin{align*}
\|\boldw\|_{\hs(B_{r_1})}^2
+\spl \sigma\boldw,\boldw \spr_{A_{r_1,r_2}}
+\spl \sigma\qcoeff\cdot\boldw,\qcoeff\cdot\boldw \spr_{A_{r_1,r_2}}
\end{align*}
is coercive due to $\sigma(\boldx)=e^{i\mu(|\boldx|)\sign\omega}$, $\mu\in[0,\mu(r_2)]$ and $\mu(r_2)\in(0,\pi)$.
Thus $T_{\hsWint}'$ is weakly coercive.
Hence there exists a compact operator $\tilde K_2\in L(\hsWint)$ such that
\begin{align*}
T_{\hsWint}:=T_{\hsWint}'+\tilde K_2
\end{align*}
is coercive and hence bijective.
Now we can define
\begin{align*}
\Tint:=\PVint-T_{\hsWint}\PWint+\PZint,
\end{align*}
which is bijective with inverse $\Tint^{-1}=\PVint-T_{\hsWint}^{-1}\PWint+\PZint$.
Subsequently let
\begin{align}\label{eq:DefTcp}
\Tcp(\boldu,v):=(\Tint\boldu,-\sigma v).
\end{align}
It follows that $T_\mathrm{cp}\in L\big(\hs_\mathrm{int}\times V\big)$ is bijective with inverse
$\Tcp^{-1}(\boldu,v)=(\Tint^{-1}\boldu,-\sigma^{-1} v)$.

\begin{theorem}\label{thm:wTc-coupled}

Let the parameters $\omega, \angvel, \sounds, \bflow, \dens, \pres, \grav$ and $\damp$ satisfy the assumptions of Section~\ref{subsec:assumptions}, in particular $\supp\bflow\subset B_{r_1}$.
Let $r_2>r_1$ be such that $\sounds, \dens \in W^{1,\infty}(B_{r_2})$, $\omega\neq0$ and $\theta$ be as defined in \eqref{eq:theta}.
Let the Hilbert spaces $\Xint$ and $\Vext$ be as defined in \eqref{eq:DefXintVext}, let the sesquilinearform $\sescp(\cdot,\cdot)$ be as defined in \eqref{eq:acp}, let $\Acp\in\blo(\Xint\times\Vext)$ be the operator associated to $\sescp(\cdot,\cdot)$ as defined in \eqref{eq:defRieszRep}, and let $\Tcp\in\blo(\Xint\times\Vext)$ be as defined in \eqref{eq:DefTcp}.
If
\begin{align*}
\|\sounds^{-1}\bflow\|_{\boldL^\infty}^2<\frac{1}{1+\tan^2\theta},
\end{align*}
then there exists $\mu$ such that $\Acp$ is weakly $\Tcp$-coercive.
\end{theorem}
\begin{proof}
For $\dislag,\dislag'\in\Xint$ we use the notation
\begin{align*}
&\compv:=\PVint\dislag, &&\compw:=\PWint\dislag, && \compz:=\PZint\dislag,\\
&\compv':=\PVint\dislag',&&\compw':=\PWint\dislag', &&\compz':=\PZint\dislag'
\end{align*}
such that $\dislag = \compv+\compw+\compz$ and $\dislag'= \compv'+\compw'+\compz'$.
Note that the space $\Xint=\hs(B_{r_2})$ differs from $\hs=\hs(\setR^3)$ and hence the subspaces $\hsVint, \hsWint, \hsZint$ differ from $\hsV, \hsW, \hsZ$ as well.
However, the role of the components $\compv, \compw, \compz$ is almost identical.
The first part of the proof is very similar to the proof of Theorem~\ref{thm:wTc-cow}.
We introduce a sufficiently small parameter $\tau\in(0,\pi/2-\theta)$ which will be specified later on,
and choose $\mu(r_2)=\pi/2-\theta-\tau\in(0,\pi/2)$.
Thus $\mu(r_2)\in(0,\pi)$, which is required for the construction of $T_{\hsWint}$.
In addition we set $\beta:=\mu_*-\pi/2+\theta+\tau$ and choose $\mu_*$ such that $\beta\in(-\pi/2,0)$ and
\begin{align}\label{eq:EstBetaCP}
0<\inf_{\boldx\in B_{r_2}^c} \inf \Re \big( i\sign\omega \, e^{-i\beta\sign\omega} \numran \, (\matii(\boldx)+i\omega\damp(\boldx)I_{3\times3})^{-1} \big).
\end{align}
Note that such a choice is possible, because $\matii(\boldx)$ is a selfadjoint matrix at each $\boldx\in B_{r_2}$ and due to $\damp\geq\dampl>0$.
We split $\Tcp^*\Acp=A_1+A_2$ into a coercive operator $A_1$ and a compact operator $A_2$.
To this end we introduce an additional small parameter $\delta>0$, which will be specified later on.
Let $\tilde K_{\sounds\sqrt{\dens}}$ be as in Theorem~\ref{thm:decomposition-Xint}.
Then similar to the proof of Theorem~\ref{thm:wTc-cow} we define $A_1$ and $A_2$ by
\begin{align*}
\spl A_1^\mathrm{int}\dislag,\dislag'\spr_{\Xint} &:=
\spl \sounds^2\div \compv,\div \compv'\spr_{B_{r_2}}
-\spl i\partial_\bflow \compv,i\partial_\bflow \compv' \spr_{B_{r_1}}
+\spl \compv,\compv'\spr_{B_{r_2}}
+\frac{1}{4\delta}\spl \tilde K_{\sounds\sqrt{\dens}}\compv,\tilde K_{\sounds\sqrt{\dens}}\compv' \spr_{\Xint}\\
&-\spl (\omega+i\partial_\bflow+i\angvel\times)\compw,i\partial_\bflow \compv' \spr_{B_{r_1}}
+\spl i\partial_\bflow \compv,(\omega+i\partial_\bflow+i\angvel\times)\compw' \spr_{B_{r_1}}
\\
&+\spl (\omega+i\partial_\bflow+i\angvel\times)\compw,(\omega+i\partial_\bflow+i\angvel\times)\compw' \spr_{B_{r_1}},\\
&+\spl (i\omega\damp+\mati)\compw,\compw'\spr_{B_{r_1}}\\
&+\spl \ol{\sigma}(i\omega\damp+\matii)\compw,\compw'\spr_{A_{r_1,r_2}}\\
&+\spl \compz,\compz' \spr_{\Xint}\\
\spl A_1\dislag,\dislag'\spr_{\Xint\times\Vext} &:=
\spl A_1^\mathrm{int}\dislag,\dislag'\spr_{\Xint}\\
&-e^{-i\mu(r_2)\sign\omega}\spl \nv\cdot\boldu,v'\spr_{H^{-1/2}(\partial B_{r_2})\times H^{1/2}(\partial B_{r_2})}\\
&-e^{-i\mu(r_2)\sign\omega}\spl v,\nv\cdot\boldu'\spr_{H^{1/2}(\partial B_{r_2})\times H^{-1/2}(\partial B_{r_2})}\\
&-\Big\spl \frac{e^{-i\mu\sign\omega}e^{2\eta}}{\dens} (\matii+i\omega\damp)^{-1} \nabla v,\nabla v'\Big\spr_{\boldL^2(B_{r_2}^c)}
+\Big\spl \frac{e^{-i\mu_*\sign\omega}e^{2\eta}}{\sounds^2\dens} v,v' \Big\spr_{\boldL^2(B_{r_2}^c)}
\end{align*}
and
\begin{align*}
\spl \op_2\dislag,\dislag'\spr_\hs &:=
-\spl \compv,\compv'\spr_{B_{r_2}}
-\frac{1}{4\delta}\spl \tilde K_{\sounds\sqrt{\dens}}\compv,\tilde K_{\sounds\sqrt{\dens}}\compv' \spr_{\Xint}
-\spl \compz,\compz' \spr_{\Xint}\\
&+\spl \div\compv,\nabla\pres\cdot\compv' \spr_{L^2(B_{r_2})}
+\spl \nabla\pres\cdot\compv,\div\compv' \spr_{L^2(B_{r_2})}
+\spl \sounds^{-2}\dens^{-1}\; \nabla\pres\cdot\compv, \nabla\pres\cdot\compv' \spr_{L^2(B_{r_2})}\\
&-\spl i\partial_\bflow \compv,(\omega+i\angvel\times)\compv' \spr_{B_{r_1}}
-\spl (\omega+i\angvel\times)\compv,i\partial_\bflow \compv' \spr_{B_{r_1}}
-\spl (i\omega\damp+\matii) \compv,\compv' \spr_{B_{r_2}}\\
&-\spl (\omega+i\partial_\bflow+i\angvel\times)\compw,(\omega+i\angvel\times) \compv' \spr_{B_{r_2}}
+\spl (\omega+i\angvel\times) \compv,(\omega+i\partial_\bflow+i\angvel\times)T'_{\hsWint}\compw' \spr_{B_{r_2}}\\
&-\spl (i\omega\damp+\mati) \compw,\compv' \spr_{B_{r_2}}
+\spl (i\omega\damp+\mati) \compv,T'_{\hsWint}\compw' \spr_{B_{r_2}}\\
&-\spl i\partial_\bflow\boldv, (\omega+i\partial_\bflow+i\angvel\times)\sigma \tilde K_1\boldw' \spr_{B_{r_1}} \\
&-\spl (\omega+i\partial_\bflow+i\angvel\times)\compw,(\omega+i\partial_\bflow+i\angvel\times)\tilde K_1\compw' \spr_{B_{r_1}}
-\spl (i\omega\damp+\mati)\compw,\tilde K_1\compw'\spr_{B_{r_1}}\\
&-\spl \ol{\sigma}(i\omega\damp+\matii)\compw,\tilde K_1\compw'\spr_{A_{r_1,r_2}}
-\spl (i\omega\damp+\matii)\compw,\check\boldv'\spr_{A_{r_1,r_2}}\\
&-\sescow(\boldv+\boldw,\tilde K_2\boldw')\\
&+\sescow^\mathrm{int}(\compv+\compw,\compz')
+\sescow^\mathrm{int}(\compz,\compv'-T_{\hsWint}\compw')
+\sescow^\mathrm{int}(\compz,\compz')\\
&+\Big\spl \frac{e^{2\eta}}{\dens} \nabla\ol{\sigma}\cdot \nabla v,v' \Big\spr_{\boldL^2(A_{r_2,r_3})}\\
&-\Big\spl \frac{e^{-i\mu_*\sign\omega}e^{2\eta}}{\sounds^2\dens} v,v' \Big\spr_{\boldL^2(A_{r_2,r_3})}
+\Big\spl \frac{e^{-i\mu\sign\omega}e^{2\eta}}{\sounds^2\dens} v,v' \Big\spr_{\boldL^2(A_{r_2,r_3})}
\end{align*}
for all $(\dislag,v),(\dislag',v')\in\Xint\times\Vext$.
Then it follows as in the proof of Theorem~\ref{thm:wTc-cow} that
\begin{align*}
\Re \left( e^{-i(\theta+\tau)\sign \omega} \spl A_1^\mathrm{int}\dislag,\dislag \spr_\hs \right)
\geq C_1 \|\dislag\|^2_{\Xint}
\end{align*}
with a constant $C_1>0$.
The coupling term, which involves both $\boldu$ and $v$, appears in
$e^{-i(\theta+\tau)\sign \omega} \spl A_1(\dislag,v),(\dislag,v) \spr_{\Xint\times\Vext}$
as
\begin{align*}
-e^{-i(\theta+\tau+\mu(r_2))\sign \omega}
2\Re\big(\spl \nv\cdot\boldu,v\spr_{H^{-1/2}(\partial B_{r_2})\times H^{1/2}(\partial B_{r_2})}\big).
\end{align*}
Due to the choice of $\mu$ it holds $-e^{-i(\theta+\tau+\mu(r_2))\sign \omega}=i\sign\omega$ and hence
\begin{align*}
\Re \left( -e^{-i(\theta+\tau+\mu(r_2))\sign \omega}
2\Re\big(\spl \nv\cdot\boldu,v\spr_{H^{-1/2}(\partial B_{r_2})\times H^{1/2}(\partial B_{r_2})}\big)\right)=0.
\end{align*}
At last we estimate the part of $A_1$ which involves only $v$.
It follows from \eqref{eq:EstBetaCP} that
\begin{align*}
\Re \left(
\Big\spl -e^{-i(\theta+\tau+\mu)\sign \omega} \frac{e^{2\eta}}{\dens} (\matii+i\omega\damp)^{-1} \nabla v,\nabla v'\Big\spr_{\boldL^2(B_{r_2}^c)}
\right)
\geq C_2 \Big\spl \frac{e^{2\eta}}{\dens} \nabla v,\nabla v'\Big\spr_{\boldL^2(B_{r_2}^c)},
\end{align*}
with a constant $C_2>0$.
Further, we note that $e^{-i(\theta+\tau+\mu_*)\sign \omega}=-i\sign\omega e^{-i\beta\sign\omega}$ and thus
$\Re\big(e^{-i(\theta+\tau+\mu_*)\sign \omega}\big)=-\sin\beta$.
Since $\beta\in(-\pi/2,0)$ we conclude
\begin{align*}
\Re \left( e^{-i(\theta+\tau)\sign \omega} \spl A_1(\dislag,v),(\dislag,v) \spr_{\Xint\times\Vext} \right)
\geq C (\|\dislag\|^2_{\Xint}+\|v\|^2_{\Vext})
\end{align*}
with $C=\min\{C_1,C_2,-\sin\beta\}>0$, i.e.\ $A_1$ is coercive.
The operator $A_2$ is compact, due to the compact embedding $\hsVint\hookrightarrow\boldL^2(B_{r_2})$ (see Theorem~\ref{thm:decomposition-Xint}, Item~\ref{it:VcompactembeddingInt}), since $\tilde K_{\sounds\sqrt{\dens}}$ is compact (see Theorem~\ref{thm:decomposition-Xint}), because $\hsZint$ is finite dimensional (see Theorem~\ref{thm:decomposition-Xint}, Item~\ref{it:ZfiniteInt}), due to the compact embedding $\boldH^1(A_{r_1,r_2})\hookrightarrow\boldL^2(A_{r_1,r_2})$ for $\check\boldv$, because $\tilde K_1$ and $\tilde K_2$ are compact, and because the embedding $\Vext\hookrightarrow L^2(A_{r_2,r_3})$ is compact.
Thus the claim is proven.
\end{proof}

\begin{proposition}\label{prop:bij-cp}
Let the assumptions of Theorem~\ref{thm:wTc-coupled} be satisfied.
Then $\Acp$ is bijective, and hence for each $\boldf\in\boldL^2_\dens(B_{r_2})$ Equation~\eqref{eq:PDE-cp} admits a unique solution $(\dislag,v)\in \Xint\times\Vext$ which depends continuously on $\boldf$.
\end{proposition}
\begin{proof}
Follows from Proposition~\ref{prop:Fredholm}, Lemma~\ref{lem:inj-cp} and Theorem~\ref{thm:wTc-coupled}.
\end{proof}


\section{Conclusion and outlook}\label{sec:conclusion}

In this article we considered the time-harmonic linear equations of stellar oscillations in $\setR^3$.
We discussed the treatment of the exterior domain and introduced a technique to couple the analysis for the interior and exterior domains.
Subsequently we established the well-posedness of Equation~\eqref{eq:PDE} in Theorem~\ref{thm:bij-full}.
In particular, the stability of \eqref{eq:sesq-a} holds without the background Equations \eqref{eq:equi} being satisfied.
In Section \ref{sec:cp} we considered the Cowling approximation with spherical symmetric parameters in the atmosphere and therein we derived a scalar equation for a potential of $\dislag$.
Subsequently we analyzed the obtained coupled system in Theorem~\ref{thm:wTc-coupled} and Proposition~\ref{prop:bij-cp}.
This new system simplifies the construction of transparent boundary conditions and leads to significant less degrees of freedom for discretizations, e.g.\ for the learned infinite elements  \cite{HohageLehrenfeldPreuss:20}.
If the factorization $(\div+\qcoeff\cdot)\dislag=e^{-\eta}\div(e^\eta\dislag)$ is applied in the entire $\setR^3$ instead of only in the atmosphere, there arise new possibilities to analyze the stellar equations, which is intended for future research.

\bibliographystyle{amsplain}
\bibliography{short_biblio}
\end{document}